\documentclass[10pt]{article}
\usepackage{amssymb,amsmath,textcomp}
\newtheorem{satz}{Theorem}[section]
\newtheorem{lemma}{Lemma}[section]

\newtheorem{korollar}{Corollary}[section]
\newtheorem{prop}{Proposition}[section]
\newcommand{\iR}{\mathbb{R}}
\newcommand{\iN}{\mathbb{N}}

\newcommand{\iC}{\mathbb{C}}

\newcommand{\oH}{\hspace*{0.39em}\raisebox{0.6ex}{\textdegree}\hspace{-0.72em}H}
\DeclareMathOperator*{\eosc}{ess\,osc}
\DeclareMathOperator*{\esup}{ess\,sup}
\DeclareMathOperator*{\einf}{ess\,inf}
\begin{document}
\begin{center}
{\bf\Large A weak Harnack inequality for fractional evolution
equations with discontinuous coefficients}
\renewcommand{\thefootnote}{\fnsymbol{footnote}}
\footnote[1]{Work partially supported by the European Community's
Human Potential Programme [Evolution Equations for Deterministic and
Stochastic Systems], contract code HPRN-CT-2002-00281, and by the
Deutsche Forschungsgemeinschaft (DFG), Bonn, Germany.}
\end{center}
\vspace{0.7em}
\begin{center}
Rico Zacher
\end{center}
{\footnotesize \noindent {\bf address:} Martin-Luther University
Halle-Wittenberg, Institute of Mathematics, Theodor-Lieser-Strasse
5, 06120 Halle, Germany, E-mail:
rico.zacher@mathematik.uni-halle.de} \vspace{0.7em}
\begin{abstract}
We study linear time fractional diffusion equations in divergence
form of time order less than one. It is merely assumed that the
coefficients are measurable and bounded, and that they satisfy a
uniform parabolicity condition. As the main result we establish for
nonnegative weak supersolutions of such problems a weak Harnack
inequality with optimal critical exponent. The proof relies on new a
priori estimates for time fractional problems and uses Moser's
iteration technique and an abstract lemma of Bombieri and Giusti,
the latter allowing to avoid the rather technically involved
approach via $BMO$. As applications of the weak Harnack inequality
we establish the strong maximum principle, continuity of weak
solutions at $t=0$, and a Liouville type theorem.
\end{abstract}
\vspace{0.7em}
\begin{center}
{\bf AMS subject classification:} 45K05, 47G20
\end{center}

\noindent{\bf Keywords:} weak Harnack inequality, Moser iterations,
fractional derivative, weak solutions, maximum principle,
subdiffusion equations, anomalous diffusion
\section{Introduction and main result}
Let $T>0$ and $\Omega$ be a bounded domain in $\iR^N$. In this paper
we are concerned with linear partial integro-differential equations
of the form
\begin{equation} \label{MProb}
\partial_t^\alpha (u-u_0)-\mbox{div}\,\big(A(t,x)Du\big)=0,\quad t\in (0,T),\,x\in
\Omega.
\end{equation}
Here $u_0=u_0(x)$ is a given initial data for $u$, $A=(a_{ij})$ is
$\iR^{N\times N}$-valued, $Du$ denotes the spatial gradient of $u$,
and $\partial_t^\alpha$ stands for the Riemann-Liouville fractional
derivation operator with respect to time of order $\alpha\in (0,1)$;
it is defined by
\[
\partial_t^\alpha v(t,x)=\partial_t\int_0^t
g_{1-\alpha}(t-\tau)v(\tau,x)\,d\tau,\quad t>0,\,x\in \Omega,
\]
where $g_\beta$ denotes the Riemann-Liouville kernel
\[
g_\beta(t)=\,\frac{t^{\beta-1}}{\Gamma(\beta)}\,,\quad
t>0,\;\beta>0.
\]

As to applications, equation (\ref{MProb}) is a special case of
problems arising in mathematical physics when describing dynamic
processes in materials with memory, e.g. in the theory of heat
conduction with memory, see \cite{JanI} and the references therein.
Time fractional diffusion equations are also used to model anomalous
diffusion, see e.g. \cite{Metz}. In this context, equations of the
type (\ref{MProb}) are termed {\em subdiffusion equations} (the time
order $\alpha$ lies in $(0,1)$; in the case $\alpha\in(1,2)$, which
is not considered here, one speaks of {\em superdiffusion
equations}. Time fractional diffusion equations of time order
$\alpha\in (0,1)$ are closely related to a class of Montroll-Weiss
continuous time random walk models where the waiting time density
behaves as $ t^{-\alpha-1}$ for $t\to \infty$, see e.g.
\cite{Hilfer1}, \cite{Hilfer2}, \cite{Metz}. Problems of the type
(\ref{MProb}) are further used to describe diffusion on fractals
(\cite{Metz}, \cite{Roma}), and they also appear in mathematical
finance, see e.g. \cite{Sca}.

Letting $\Omega_T=(0,T)\times \Omega$ we will assume that
\begin{itemize}
\item [{\bf (H1)}] $A\in L_\infty(\Omega_T;\iR^{N\times
N})$, and
\[
\sum_{i,j=1}^N|a_{ij}(t,x)|^2\le \Lambda^2,\quad \mbox{for
a.a.}\;(t,x)\in \Omega_T.
\]
\item [{\bf (H2)}] There exists $\nu>0$ such that
\[
\big(A(t,x)\xi|\xi\big)\ge \nu|\xi|^2,\quad\mbox{for a.a.}\;
(t,x)\in\Omega_T,\; \mbox{and all}\;\xi\in \iR^N.
\]
\item [{\bf (H3)}] $u_0\in L_2(\Omega)$.
\end{itemize}
We say that a function $u$ is a {\em weak solution (subsolution,
supersolution)} of (\ref{MProb}) in $\Omega_T$, if $u$ belongs to
the space
\begin{align*}
Z_\alpha:=\{&\,v\in
L_{\frac{2}{1-\alpha},\,w}([0,T];L_2(\Omega))\cap
L_2([0,T];H^1_2(\Omega))\;
\mbox{such that}\;\\
&\;\;g_{1-\alpha}\ast v\in C([0,T];L_2(\Omega)),
\;\mbox{and}\;(g_{1-\alpha}\ast v)|_{t=0}=0\},
\end{align*}
and for any nonnegative test function
\[
\eta\in \oH^{1,1}_2(\Omega_T):=H^1_2([0,T];L_2(\Omega))\cap
L_2([0,T];\oH^1_2(\Omega)) \quad\quad
\Big(\oH^1_2(\Omega):=\overline{C_0^\infty(\Omega)}\,{}^{H^1_2(\Omega)}\Big)
\]
with $\eta|_{t=T}=0$ there holds
\begin{equation} \label{BWF}
\int_{0}^{T} \int_\Omega \Big(-\eta_t [g_{1-\alpha}\ast (u-u_0)]+
(ADu|D \eta)\Big)\,dxdt= \,(\le,\,\ge )\,0.
\end{equation}
Here $L_{p,\,w}$ denotes the weak $L_p$ space and $f_1\ast f_2$ the
convolution on the positive halfline with respect to time, that is
$(f_1\ast f_2)(t)=\int_0^t f_1(t-\tau)f_2(\tau)\,d\tau$, $t\ge 0$.

Weak solutions of (\ref{MProb}) in the class $Z_\alpha$ have been
constructed in \cite{ZWH}. Notice also that the function $u_0$ plays
the role of the initial data for $u$, at least in a weak sense. In
case of sufficiently smooth functions $u$ and
$g_{1-\alpha}\ast(u-u_0)$ the condition $(g_{1-\alpha}\ast
u)|_{t=0}=0$ implies $u|_{t=0}=u_0$, see \cite{ZWH}.

To formulate our main result, let $B(x,r)$ denote the open ball with
radius $r>0$ centered at $x\in \iR^N$. By $\mu_N$ we mean the
Lebesgue measure in $\iR^N$. For $\delta\in(0,1)$, $t_0\ge 0$,
$\tau>0$, and a ball $B(x_0,r)$, define the boxes
\begin{align*}
Q_-(t_0,x_0,r)&=(t_0,t_0+\delta\tau r^{2/\alpha})\times B(x_0,\delta r),\\
Q_+(t_0,x_0,r)&=(t_0+(2-\delta)\tau r^{2/\alpha},t_0+2\tau
r^{2/\alpha})\times B(x_0,\delta r).
\end{align*}
\begin{satz} \label{localweakHarnack}
Let $\alpha\in(0,1)$, $T>0$, and $\Omega\subset \iR^N$ be a
bounded domain. Suppose the assumptions (H1)--(H3) are satisfied.
Let further $\delta\in(0,1)$, $\eta>1$, and $\tau>0$ be fixed.
Then for any $t_0\ge 0$ and $r>0$ with $t_0+2\tau r^{2/\alpha}\le
T$, any ball $B(x_0,\eta r)\subset\Omega$, any
$0<p<\frac{2+N\alpha}{2+N\alpha-2\alpha}$, and any nonnegative
weak supersolution $u$ of (\ref{MProb}) in $(0,t_0+2\tau
r^{2/\alpha})\times B(x_0,\eta r)$ with $u_0\ge 0$ in $B(x_0,\eta
r)$, there holds
\begin{equation} \label{localwHarnackF}
\Big(\frac{1}{\mu_{N+1}\big(Q_-(t_0,x_0,r)\big)}\,\int_{Q_-(t_0,x_0,r)}u^p\,d\mu_{N+1}\Big)^{1/p}
\le C \einf_{Q_+(t_0,x_0,r)} u,
\end{equation}
where the constant $C=C(\nu,\Lambda,\delta,\tau,\eta,\alpha,N,p)$.
\end{satz}
Theorem \ref{localweakHarnack} states that nonnegative weak
supersolutions of (\ref{MProb}) satisfy a weak form of Harnack
inequality in the sense that we do not have an estimate for the
supremum of $u$ on $Q_-(t_0,x_0,r)$ but only an $L_p$ estimate. We
also show that the critical exponent
$\frac{2+N\alpha}{2+N\alpha-2\alpha}$ is optimal, i.e. the
inequality fails to hold for $p\ge
\frac{2+N\alpha}{2+N\alpha-2\alpha}$.

Theorem \ref{localweakHarnack} can be viewed as the time fractional
analogue of the corresponding result in the classical parabolic case
$\alpha=1$, see e.g. \cite[Theorem 6.18]{Lm} and \cite{Trud}.
Sending $\alpha\to 1$ in the expression for the critical exponent
yields $1+2/N$, which is the well-known critical exponent for the
heat equation. We would like to point out that the statement of
Theorem \ref{localweakHarnack} remains valid for (appropriately
defined) weak supersolutions of (\ref{MProb}) on $(t_0,t_0+2\tau
r^{2/\alpha})\times B(x_0,\eta r)$ which are nonnegative on
$(0,t_0+2\tau r^{2/\alpha})\times B(x_0,\eta r)$. Here the global
positivity assumption cannot be replaced by a local one, as simple
examples show, cf. \cite{Za3}. This significant difference to the
case $\alpha=1$ is due to the non-local nature of
$\partial_t^\alpha$. The same phenomenon is known for
integro-differential operators like $(-\Delta)^\alpha$ with
$\alpha\in (0,1)$, see e.g. \cite{Kass1}.

As a simple consequence of the weak Harnack inequality we derive the
strong maximum principle for weak subsolutions of (\ref{MProb}), see
Theorem \ref{strongmax} below. The weak maximum principle has been
proven in \cite{Za2}, even in a more general setting.

In the classical parabolic case boundedness and the weak (or full)
Harnack inequality imply an H\"older estimate for weak solutions,
cf.\ \cite{DB}, \cite{LSU}, \cite{Lm}, \cite{Moser64}. We also refer
to \cite{GilTrud} and \cite{Mosell} for the elliptic case. In the
present situation one cannot argue anymore as in the classical
parabolic case, due to the global positivity assumption in Theorem
\ref{localweakHarnack}. The same problem arises for the fractional
Laplacian, see \cite{Silv}. However, in our case it is possible to
establish at least continuity at $t=0$. This is done in Theorem
\ref{Hoeldert=0} in the case $u_0=0$. It is shown that in this case
any bounded weak solution $u$ of (\ref{MProb}) is continuous at
$(0,x_0)$ for all $x_0\in \Omega$ and $\lim_{(t,x)\to
(0,x_0)}u(t,x)=0$. Thus for such weak solutions the initial
condition $u|_{t=0}=0$ is satisfied in the classical sense.

As a further consequence of the weak Harnack inequality we obtain a
theorem of Liouville type, see Corollary \ref{Liouville} below. It
states that any bounded weak solution of (\ref{MProb}) on
$\iR_+\times \iR^N$ with $u_0=0$ vanishes a.e. on $\iR_+\times
\iR^N$.

The proof of Theorem \ref{localweakHarnack} relies on new a priori
estimates for time fractional problems, which are derived by means
of the fundamental identity (\ref{fundidentity}) (see below) for the
regularized fractional derivative. It further uses Moser's iteration
technique and an elementary but subtle lemma of Bombieri and Giusti
\cite{BomGiu}, which allows to avoid the rather technically involved
approach via $BMO$-functions. This simplification is already of
great significance in the classical parabolic case, see Moser
\cite{Moser71} and Saloff-Coste \cite{SalCoste}.

One of the technical difficulties in deriving the desired estimates
in the weak setting is to find an appropriate time regularization of
the problem. In the case $\alpha=1$ this can be achieved by means of
Steklov averages in time. In the time fractional case this method
does not work anymore, since Steklov average operators and
convolution do not commute. It turns out that instead one can use
the Yosida approximation of the fractional derivative, which leads
to a regularization of the kernel $g_{1-\alpha}$. This method has
already been used in \cite{Grip1}, \cite{VZ}, \cite{ZWH}, and
\cite{Za2}.

We point out that the results obtained in this paper can be easily
generalized to quasilinear equations of the form
\begin{equation} \label{QProb}
\partial_t^\alpha (u-u_0)-\mbox{div}\,a(t,x,u,Du)=b(t,x,u,Du),\;t\in (0,T),\,x\in \Omega,
\end{equation}
with suitable structure conditions on the functions $a$ and $b$.
This is possible, as also known from the elliptic and the classical
parabolic case, since the test function method used in the proof of
Theorem \ref{localweakHarnack} does not depend so much on the
linearity of the differential operator w.r.t.\ the spatial variables
but on a certain nonlinear structure, cf. \cite{GilTrud}, \cite{Lm},
\cite{Trud}, and \cite{Za2}.

In the literature there exist many papers where equations of the
type (\ref{MProb}), as well as nonlinear or abstract variants of
them are studied in a {\em strong} setting, assuming more smoothness
on the coefficients and nonlinearities, see e.g. \cite{Ba},
\cite{CLS}, \cite{CP2}, \cite{Koch}, \cite{Grip1}, \cite{JanI},
\cite{ZEQ}, \cite{ZQ}. Concerning the {\em weak} setting described
above one finds only a few results. Existence of weak solutions has
been shown in \cite{ZWH} in an abstract setting for a more general
class of kernels. Boundedness of weak solutions has been obtained in
\cite{Za2} in the quasilinear case by means of the De Giorgi
technique. With the weak Harnack inequality, the present paper
establishes a key result towards a De Giorgi-Nash-Moser theory for
time fractional evolution equations in divergence form of order
$\alpha\in(0,1)$.

We further remark that in the purely time-dependent case, that is
for scalar equations of the form
\[
\partial_t^\alpha (u-u_0)+\sigma u=0,\quad t\in (0,T),
\]
with $\sigma\ge 0$, a weak Harnack inequality with optimal exponent
$1/(1-\alpha)$ has been proven in \cite{Za} for nonnegative
supersolutions. Recently, the full Harnack inequality for
nonnegative solutions has been established in \cite{Za3}. This,
together with the above results, indicates that the full Harnack
inequality should also hold for nonnegative solutions of
(\ref{MProb}), which is still an open problem, even in the case
$A(t,x)\equiv Id$.

The paper is organized as follows. In Section 2 we collect the basic
tools needed for the proof of Theorem \ref{localweakHarnack}. These
include two abstract lemmas on Moser iterations and the lemma of
Bombieri and Giusti. We further explain the approximation method for
the fractional derivation operator and state the fundamental
identity (\ref{fundidentity}), which is frequently used in Section
3, where we give the proof of the main result. In Section 4 we show
that the critical exponent in Theorem \ref{localweakHarnack} is
optimal. Finally, Section 5 is devoted to applications of the weak
Harnack inequality.
\section{Preliminaries}
\subsection{Moser iterations and an abstract lemma of Bombieri and Giusti}
Throughout this subsection $U_\sigma$, $0<\sigma\le 1$, will denote
a collection of measurable subsets of a fixed finite measure space
endowed with a measure $\mu$, such that $U_{\sigma'}\subset
U_\sigma$ if $\sigma'\le \sigma$. For $p\in (0,\infty)$ and
$0<\sigma\le 1$, $L_p(U_\sigma)$ stands for the Lebesgue space
$L_p(U_\sigma,d\mu)$ of all $\mu$-measurable functions
$f:U_\sigma\rightarrow \iR$ with
$|f|_{L_p(U_\sigma)}:=(\int_{U_\sigma}|f|^p\,d\mu)^{1/p}<\infty$.

The following two lemmas are basic to Moser's iteration technique.
The arguments in their proofs have been repeatedly used in the
literature (see e.g. \cite{GilTrud}, \cite{Lm}, \cite{Mosell},
\cite{Moser64}, \cite{SalCoste}, \cite{Trud}), so it is worthwhile
to formulate them as lemmas in abstract form, also for future
reference. We provide proofs for the sake of completeness.

The first Moser iteration result reads as follows, see also
\cite[Lemma 2.3]{CZ}.
\begin{lemma} \label{moserit1}
Let $\kappa>1$, $\bar{p}\ge 1$, $C\ge 1$, and $\gamma>0$. Suppose
$f$ is a $\mu$-measurable function on $U_1$ such that
\begin{equation} \label{mositer1}
|f|_{L_{\beta\kappa}(U_{\sigma'})}\le
\Big(\frac{C(1+\beta)^{\gamma}}{(\sigma-\sigma')^{\gamma}}\Big)^{1/\beta}\,|f|_{L_{\beta}(U_{\sigma})},
\quad 0<\sigma'<\sigma\le 1,\;\beta>0.
\end{equation}
Then there exist constants $M=M(C,\gamma,\kappa,\bar{p})$ and
$\gamma_0=\gamma_0(\gamma,\kappa)$ such that
\[
\esup_{U_{\delta}}{|f|} \le
\Big(\frac{M}{(1-\delta)^{\gamma_0}}\Big)^{1/p}
|f|_{L_{p}(U_1)}\quad \mbox{for all}\;\;\delta\in(0,1),\;p\in
(0,\bar{p}] .\]
\end{lemma}
{\em Proof:} For $q>0$ and $0<\sigma\le 1$, let
\[ \Phi(q,\sigma)=(\int_{U_\sigma} |f|^q\,d\mu)^{1/q}.\]
Let $0<p\le \bar{p}$ and $\delta\in (0,1)$. Set $p_i=p\kappa^{i}$,
$i=0,1,\ldots$ and define the sequence $\{\sigma_i\}$,
$i=0,1,\ldots$, by $\sigma_0=1$ and $\sigma_i=1-\sum_{j=1}^i 2^{-j}
(1-\delta)$, $i=1,2,\ldots$; observe that
$1=\sigma_0>\sigma_1>\ldots>\sigma_i>\sigma_{i+1}>\delta$ as well as
$\sigma_{i-1}-\sigma_{i}=2^{-i}(1-\delta)$, $i\ge 1$. Suppose now
$n\in \iN$. By using (\ref{mositer1}) with $\beta=p_i$,
$i=0,1,\ldots,n-1$, we obtain
\begin{align*}
\Phi(p_n,\delta)& \;\le
\Phi(p_n,\sigma_n)\;=\;\Phi(p_{n-1}\kappa,\sigma_n)\;\le\;
\Big(\frac{C(1+p\kappa^{n-1})^{\gamma}}{[2^{-n}(1-\delta)]^{\gamma}}\Big)^
{\frac{1}{p}\,\kappa^{-(n-1)}}\Phi(p_{n-1},\sigma_{n-1})\\ & \;\le
\Big(\frac{C(2\bar{p}\kappa^{n-1})^{\gamma}}{[2^{-n}(1-\delta)]^{\gamma}}\Big)^
{\frac{1}{p}\,\kappa^{-(n-1)}}\Phi(p_{n-1},\sigma_{n-1})\\ &\;\le
\Big(\frac{\tilde{C}(C,\bar{p},\gamma)^n
\kappa^{\gamma(n-1)}}{(1-\delta)^{\gamma}}\Big)^
{\frac{1}{p}\,\kappa^{-(n-1)}}\Phi(p_{n-1},\sigma_{n-1})\;\le\;\ldots\\
& \;\le \Big(\tilde{C}^{\sum_{j=0}^{n-1}
(j+1)\kappa^{-j}}\kappa^{\gamma\sum_{j=0}^{n-1}
j\kappa^{-j}}(1-\delta)^{-\gamma\sum_{j=0}^{n-1}
\kappa^{-j}}\Big)^{1/p}\,\Phi(p_0,\sigma_0)\\ & \; \le
\Big(\frac{M(C,\bar{p},\gamma,\kappa)}{(1-\delta)^{\frac{\gamma\kappa}{\kappa-1}}}\Big)^
{1/p}\,\Phi(p,1).
\end{align*}
We let now $n$ tend to $\infty$ and use the fact that
\[ \lim_{n\to\infty}\Phi(p_n,\delta)=\esup_{U_{\delta}}{|f|}\]
to get
\[ \esup_{U_{\delta}}{|f|} \le
\Big(\frac{M(C,\bar{p},\gamma,\kappa)}{(1-\delta)^{\frac{\gamma\kappa}{\kappa-1}}}\Big)^
{1/p}\,|f|_{L_{p}(U_1)}.\] Hence the proof is complete. $\square$

$\mbox{}$

\noindent The second Moser iteration result is the following, see
also \cite[Lemma 2.5]{CZ}.
\begin{lemma} \label{moserit2}
Assume that $\mu(U_1)\le 1$. Let $\kappa>1$, $0<p_0<\kappa$, and
$C\ge 1,\,\gamma>0$. Suppose $f$ is a $\mu$-measurable function on
$U_1$ such that
\begin{equation} \label{mositer2}
|f|_{L_{\beta\kappa}(U_{\sigma'})}\le
\Big(\frac{C}{(\sigma-\sigma')^{\gamma}}\Big)^{1/\beta}\,|f|_{L_{\beta}(U_{\sigma})},
\quad 0<\sigma'<\sigma\le 1,\;0<\beta\le \frac{p_0}{\kappa}<1.
\end{equation}
Then there exist constants $M=M(C,\gamma,\kappa)$ and
$\gamma_0=\gamma_0(\gamma,\kappa)$ such that
\[ |f|_{L_{p_0}(U_{\delta})}\le
\Big(\frac{M}{(1-\delta)^{\gamma_0}}\Big)^{1/p-1/p_0}
|f|_{L_{p}(U_1)}\quad \mbox{for all}\;\;\delta\in(0,1),\;p\in
(0,\frac{p_0}{\kappa}] .\]
\end{lemma}
{\em Proof:} Set $p_i=p_0 \kappa^{-i}$, $i=1,2,\ldots$. Given
$\delta\in(0,1)$ we take again the sequence $\{\sigma_i\}$,
$i=0,1,2,\ldots$, defined by $\sigma_0=1$ and
$\sigma_i=1-\sum_{j=1}^i 2^{-j} (1-\delta)$, $i\ge 1$. Suppose now
$n\in \iN$. By using (\ref{mositer2}) with $\beta=p_i$,
$i=1,\ldots,n$, we obtain
\begin{align*}
\Phi(p_0,\delta) & \le \; \Phi(p_0,\sigma_n)\;=\;\Phi(p_1
\kappa,\sigma_n)\;\le
\;\frac{C^{\kappa/p_0}}{[2^{-n}(1-\delta)]^{\gamma
 \kappa/p_0}}\;\Phi(p_1,\sigma_{n-1})\\
& \le \; \frac{C^{\kappa/p_0}}{[2^{-n}(1-\delta)]^{\gamma
 \kappa/p_0}}\; \frac{C^{\kappa^2/p_0}}{[2^{-(n-1)}(1-\delta)]^{\gamma
 \kappa^2/p_0}}\;\Phi(p_2,\sigma_{n-2})\;\le\;\dots\\
& \le \; \frac{C^{\frac{1}{p_0}\,(\kappa+\kappa^2+\ldots+\kappa^n)}}
{2^{-\frac{\gamma}{p_0}\,(n\kappa+(n-1)\kappa^2+\ldots+2\kappa^{n-1}+\kappa^n)}
(1-\delta)^{\frac{\gamma}{p_0}\,(\kappa+\kappa^2+\ldots+\kappa^n)}}\;\Phi(p_n,\sigma_0).
\end{align*}
Since $p_i=p_0 \kappa^{-i}$, we have
\[ \frac{1}{p_0}\,\sum_{j=1}^n \kappa^j = \frac{\kappa(\kappa^n-1)}{p_0(\kappa-1)}
= \frac{\kappa}{p_0(\kappa-1)}\;(\frac{p_0}{p_n}-1)=
\frac{\kappa}{\kappa-1}\;(\frac{1}{p_n}-\frac{1}{p_0}).\] Employing
the formula
\[ \sum_{j=1}^n j
\kappa^{j-1}=\frac{1-(n+1)\kappa^n+n\kappa^{n+1}}{(\kappa-1)^2}\] we
have further
\begin{align*}
\sum_{j=1}^n (n+1-j)\kappa^j & = \; (n+1)\sum_{j=1}^n
\kappa^j-\sum_{j=1}^n j \kappa^j\\ & =
\;(n+1)\kappa\,\frac{\kappa^n-1}{\kappa-1}-\kappa\;
\frac{1-(n+1)\kappa^n+n\kappa^{n+1}}{(\kappa-1)^2}\\ & = \;
\kappa\;\frac{\kappa^{n+1}-(n+1)\kappa+n}{(\kappa-1)^2}\;\le \;
\frac{\kappa}{(\kappa-1)^2}\kappa^{n+1}\\ & \le \;
\frac{\kappa^3}{(\kappa-1)^3}\;(\kappa^n-1)\;\le\;\frac{\kappa^3}{(\kappa-1)^3}\;
(\frac{p_0}{p_n}-1),
\end{align*}
which yields
\[ \frac{1}{p_0}\,\sum_{j=1}^n (n+1-j)\kappa^j \le \frac{\kappa^3}{(\kappa-1)^3}\;
(\frac{1}{p_n}-\frac{1}{p_0}).\] Therefore
\[ \Phi(p_0,\delta) \le \Big[
\frac{2^{\frac{\gamma
\kappa^3}{(\kappa-1)^3}}C^{\frac{\kappa}{\kappa-1}}}
{(1-\delta)^{\frac{\gamma\kappa}{\kappa-1}}}\Big]^
{\frac{1}{p_n}-\frac{1}{p_0}} \Phi(p_n,\sigma_0).
\]
Given $p\in(0,p_0/\kappa]$ there exists $n\ge 2$ such that $p_n<p\le
p_{n-1}$. We then have
\begin{align*}
\frac{1}{p_n}-\frac{1}{p_0} & \; =\frac{\kappa^n-1}{p_0}\;\le
\;\frac{\kappa^n+\kappa^{n-1}-\kappa-1}{p_0}\;=\;
\frac{(1+\kappa)(\kappa^{n-1}-1)}{p_0}\\ & \; =
(1+\kappa)(\frac{1}{p_{n-1}}-\frac{1}{p_0})\;\le\;(1+\kappa)(\frac{1}{p}-\frac{1}{p_0}),
\end{align*}
as well as
\[ \Phi(p_n,\sigma_0)=\Phi(p_n,1)\le \Phi(p,1),\]
by H\"older's inequality and the assumption $\mu(U_1)\le 1$. All in
all, we obtain
\[ \Phi(p_0,\delta)\le \Big[
\frac{2^{\frac{\gamma
\kappa^3}{(\kappa-1)^3}}C^{\frac{\kappa}{\kappa-1}}}
{(1-\delta)^{\frac{\gamma\kappa}{\kappa-1}}}\Big]^
{(1+\kappa)(\frac{1}{p}-\frac{1}{p_0})}\Phi(p,1),\] which proves the
lemma. $\square$

$\mbox{}$

The following abstract lemma is due to Bombieri and Giusti
\cite{BomGiu}. For a proof we also refer to \cite[Lemma
2.2.6]{SalCoste} and \cite[Lemma 2.6]{CZ}
\begin{lemma} \label{abslemma}
Let $\delta,\,\eta\in(0,1)$, and let $\gamma,\,C$ be positive
constants and $0<\beta_0\le \infty$. Suppose $f$ is a positive
$\mu$-measurable function on $U_1$ which satisfies the following two
conditions:

(i)
\[
|f|_{L_{\beta_0}(U_{\sigma'})}\le
[C(\sigma-\sigma')^{-\gamma}\mu(U_1)^{-1}]^{1/\beta-1/\beta_0}|f|_{L_{\beta}(U_{\sigma})},
\]
for all $\sigma,\,\sigma',\,\beta$ such that $0<\delta\le
\sigma'<\sigma\le 1$ and $0<\beta\le \min\{1,\eta\beta_0\}$.

(ii)
\[
\mu(\{\log f>\lambda\})\le C\mu(U_1)\lambda^{-1}
\]
for all $\lambda>0$.

Then
\[
|f|_{L_{\beta_0}(U_{\delta})}\le M \mu(U_1)^{1/\beta_0},
\]
where $M$ depends only on $\delta,\,\eta,\,\gamma,\,C$, and
$\beta_0$.
\end{lemma}
\subsection{The Yoshida approximation of the fractional
derivation operator} \label{SecYos}
Let $0<\alpha<1$, $1\le p<\infty$, $T>0$, and $X$ be a real Banach
space. Then the fractional derivation operator $B$ defined by
\[ B u=\,\frac{d}{dt}\,(g_{1-\alpha}\ast u),\;\;D(B)=\{u\in L_p([0,T];X):\,g_{1-\alpha}\ast u\in \mbox{}_0 H^1_p([0,T];X)\},
\]
where the zero means vanishing at $t=0$, is known to be
$m$-accretive in $L_p([0,T];X)$, cf. \cite{Phil1}, \cite{CP}, and
\cite{Grip1}. Its Yosida approximations $B_{n}$, defined by
$B_{n}=nB(n+B)^{-1},\,n\in \iN$, enjoy the property that for any
$u\in D(B)$, one has $B_{n}u\rightarrow Bu$ in $L_p([0,T];X)$ as
$n\to \infty$. Further, one has the representation
\begin{equation} \label{Yos}
B_n u=\,\frac{d}{dt}\,(g_{1-\alpha,n}\ast u),\quad u\in
L_p([0,T];X),\;n\in \iN,
\end{equation}
where $g_{1-\alpha,n}=n s_{\alpha,n}$, and $s_{\alpha,n}$ is the
unique solution of the scalar-valued Volterra equation
\[
s_{\alpha,n}(t)+n(s_{\alpha,n}\ast g_\alpha)(t)=1,\quad
t>0,\;n\in\iN,
\]
see e.g. \cite{VZ}. Let $h_{\alpha,n}\in L_{1,\,loc}(\iR_+)$ be the
resolvent kernel associated with $ng_\alpha$, that is
\begin{equation} \label{hndef}
h_{\alpha,n}(t)+n(h_{\alpha,n}\ast g_\alpha)(t)=ng_\alpha(t),\quad
t>0,\;n\in\iN.
\end{equation}
Convolving (\ref{hndef}) with $g_{1-\alpha}$ and using $g_\alpha\ast
g_{1-\alpha}=1$, we obtain
\[
(g_{1-\alpha}\ast h_{\alpha,n})(t)+n([g_{1-\alpha}\ast
h_{\alpha,n}]\ast g_\alpha)(t)=n,\quad t>0,\;n\in\iN.
\]
Hence
\begin{equation} \label{gnprop}
g_{1-\alpha,n}=ns_{\alpha,n}=g_{1-\alpha}\ast h_{\alpha,n},\quad
n\in \iN.
\end{equation}
The kernels $g_{1-\alpha,n}$ are nonnegative and nonincreasing for
all $n\in\iN$, and they belong to $H^1_1([0,T])$, cf. \cite{JanI}
and \cite{VZ}. Note that for any function $f\in L_p([0,T];X)$, $1\le
p<\infty$, there holds $h_{\alpha,n}\ast f\to f$ in $L_p([0,T];X)$
as $n\to \infty$. In fact, setting $u=g_\alpha\ast f$, we have $u\in
D(B)$, and
\[
B_n u=\,\frac{d}{dt}\,(g_{1-\alpha,n}\ast
u)=\,\frac{d}{dt}\,(g_{1-\alpha}\ast g_\alpha\ast h_{\alpha,n}\ast
f)=h_{\alpha,n}\ast f\,\to\,Bu=f\quad\mbox{in}\;L_p([0,T];X)
\]
as $n\to \infty$. In particular, $g_{1-\alpha,n}\to g_{1-\alpha}$ in
$L_1([0,T])$ as $n\to \infty$.

We next state a fundamental identity for integro-differential
operators of the form $\frac{d}{dt}(k\ast u)$, cf.\ also \cite{Za2}.
Suppose $k\in H^1_1([0,T])$ and $H\in C^1(\iR)$. Then it follows
from a straightforward computation that for any sufficiently smooth
function $u$ on $(0,T)$ one has for a.a. $t\in (0,T)$,
\begin{align} \label{fundidentity}
H'(u(t))&\frac{d}{dt}\,(k \ast u)(t) =\;\frac{d}{dt}\,(k\ast
H(u))(t)+
\Big(-H(u(t))+H'(u(t))u(t)\Big)k(t) \nonumber\\
 & +\int_0^t
\Big(H(u(t-s))-H(u(t))-H'(u(t))[u(t-s)-u(t)]\Big)[-\dot{k}(s)]\,ds,
\end{align}
where $\dot{k}$ denotes the derivative of $k$. In particular this
identity applies to the Yosida approximations of the fractional
derivation operator. We remark that an integrated version of
(\ref{fundidentity}) can be found in \cite[Lemma 18.4.1]{GLS}.
Observe that the last term in (\ref{fundidentity}) is nonnegative in
case $H$ is convex and $k$ is nonincreasing.

The subsequent two lemmas are also obtained by simple algebra.
\begin{lemma} \label{comm}
Let $T>0$ and $\alpha\in (0,1)$. Suppose that $v\in
{}_0H^1_1([0,T])$ and $\varphi\in C^1([0,T])$. Then
\[
\big(g_\alpha\ast(\varphi \dot{v}))(t)=\varphi(t)(g_\alpha\ast
\dot{v})(t)+\int_0^t
v(\sigma)\partial_\sigma\big(g_\alpha(t-\sigma)
[\varphi(t)-\varphi(\sigma)]\big)\,d\sigma,\;\;\mbox{a.a.}\;t\in
(0,T).
\]
If in addition $v$ is nonnegative and $\varphi$ is nondecreasing
there holds
\[
\big(g_\alpha\ast(\varphi \dot{v}))(t)\ge \varphi(t)(g_\alpha\ast
\dot{v})(t)-\int_0^t g_\alpha(t-\sigma)
\dot{\varphi}(\sigma)v(\sigma)\,d\sigma,\;\;\mbox{a.a.}\;t\in
(0,T).
\]
\end{lemma}
\begin{lemma} \label{comm2}
Let $T>0$, $k\in H^1_1([0,T])$, $v\in L_1([0,T])$, and $\varphi\in
C^1([0,T])$. Then
\[
\varphi(t)\,\frac{d}{dt}\,(k\ast v)(t)=\,\frac{d}{dt}\,\big(k\ast
[\varphi v]\big)(t)+\int_0^t
\dot{k}(t-\tau)\big(\varphi(t)-\varphi(\tau)\big)v(\tau)\,d\tau,\;\;\mbox{a.a.}\;t\in
(0,T).
\]

\end{lemma}
\subsection{An embedding result and a weighted Poincar\'e inequality}
Let $T>0$ and $\Omega$ be a bounded domain in $\iR^N$. For $1<p\le
\infty$ we define the space
\begin{equation} \label{Vdef}
V_p:=V_p([0,T]\times \Omega)=L_{2p}([0,T];L_2(\Omega))\cap
L_2([0,T];H^1_2(\Omega)),
\end{equation}
endowed with the norm
\[
|u|_{V_p([0,T]\times \Omega)}:=|u|_{L_{2p}([0,T];L_2(\Omega))}
+|Du|_{L_2([0,T];L_2(\Omega))}.
\]
Set
\begin{equation} \label{kappa}
\kappa:=\kappa_p:=\,\frac{2p+N(p-1)}{2+N(p-1)}
\end{equation}
with $\kappa_\infty=1+2/N$. Then $V_p\hookrightarrow
L_{2\kappa}([0,T]\times\Omega)$, and
\begin{equation} \label{Vembedding}
|u|_{L_{2\kappa}([0,T]\times\Omega)}\le C(N,p)|u|_{V_p([0,T]\times
\Omega)},
\end{equation}
for all $u\in V_p\cap L_2([0,T];\oH^1_2(\Omega))$. This is a
consequence of the Gagliardo-Nirenberg and H\"older's inequality.
The case $p=\infty$ is contained, e.g., in \cite[p.\ 74 and
75]{LSU}. The proof given there easily extends to the general case.
For a more general embedding result (without proof) we also refer to
\cite[Section 2]{Za2}.

The following result can be found in \cite[Lemma 3]{Moser64}, see
also \cite[Lemma 6.12]{Lm}.
\begin{prop} \label{WeiPI}
Let $\varphi\in C(\iR^N)$ with non-empty compact support of diameter
$d$ and assume that $0\le \varphi\le 1$. Suppose that the domains
$\{x\in\iR^N:\varphi(x)\ge a\}$ are convex for all $a\le 1$. Then
for any function $u\in H^{1}_2(\iR^N)$,
\[
\int_{\iR^N} \big(u(x)-u_\varphi\big)^2 \varphi(x)\,dx \le \,\frac{2
d^2\mu_N(\mbox{{\em supp}}\,\varphi)}{|\varphi|_{L_1(\iR^N)}}\,
\int_{\iR^N} |Du(x)|^2 \varphi(x)\,dx,
\]
where
\[
u_\varphi=\frac{\int_{\iR^N} u(x)\varphi(x)\,dx}{\int_{\iR^N}
\varphi(x)\,dx}.
\]
\end{prop}
\section{Proof of the main result}
\subsection{The regularized weak formulation, time shifts, and
scalings}\label{SSS}
The following lemma is basic to deriving {\em a priori} estimates
for weak (sub-/super-) solutions of (\ref{MProb}). It provides an
equivalent weak formulation of (\ref{MProb}) where the singular
kernel $g_{1-\alpha}$ is replaced by the more regular kernel
$g_{1-\alpha,n}$ ($n\in\iN$) given in (\ref{gnprop}). In what
follows the kernels $h_n:=h_{\alpha,n}$, $n\in\iN$, are defined as
in Section \ref{SecYos}.
\begin{lemma} \label{LemmaReg}
Let $\alpha\in (0,1)$, $T>0$, and $\Omega\subset \iR^N$ be a bounded
domain. Suppose the assumptions (H1)--(H3) are satisfied. Then $u\in
Z_\alpha$ is a weak solution (subsolution, supersolution) of
(\ref{MProb}) in $\Omega_T$ if and only if for any nonnegative
function $\psi\in \oH^1_2(\Omega)$ one has
\begin{equation} \label{LemmaRegF}
\int_\Omega \Big(\psi \partial_t[g_{1-\alpha,n}\ast
(u-u_0)]+(h_n\ast [ADu]|D\psi)\Big)\,dx\nonumber\\
=\,(\le,\,\ge)\,0,\quad\mbox{a.a.}\;t\in (0,T),\,n\in \iN.
\end{equation}
\end{lemma}
For a proof we refer to Lemma 3.1 in \cite{Za2}, where a more
general situation is considered with a slightly different function
space for the solution. The proof of Lemma \ref{LemmaReg} is
analogous.

Let $u\in Z_\alpha$ be a weak supersolution of (\ref{MProb}) in
$\Omega_T$ and assume that $u_0\ge 0$ in $\Omega$. Then Lemma
\ref{LemmaReg} and positivity of $g_{1-\alpha,n}$ imply that
\begin{equation} \label{u0weg}
\int_\Omega \Big(\psi \partial_t(g_{1-\alpha,n}\ast u)+(h_n\ast
[ADu]|D\psi)\Big)\,dx \ge \,0,\quad\mbox{a.a.}\;t\in (0,T),\,n\in
\iN,
\end{equation}
for any nonnegative function $\psi\in \oH^1_2(\Omega)$.

Let now $t_1\in (0,T)$ be fixed. For $t\in (t_1,T)$ we introduce the
shifted time $s=t-t_1$ and set $\tilde{f}(s)=f(s+t_1)$, $s\in
(0,T-t_1)$, for functions $f$ defined on $(t_1,T)$. From the
decomposition
\[
(g_{1-\alpha,n}\ast u)(t,x)=\int_{t_1}^t
g_{1-\alpha,n}(t-\tau)u(\tau,x)\,d\tau+\int_{0}^{t_1}
g_{1-\alpha,n}(t-\tau)u(\tau,x)\,d\tau,\quad t\in (t_1,T),
\]
we then deduce that
\begin{equation} \label{shiftprop}
\partial_t(g_{1-\alpha,n}\ast u)(t,x)=\partial_s(g_{1-\alpha,n}\ast
\tilde{u})(s,x)+\int_0^{t_1}\dot{g}_{1-\alpha,n}(s+t_1-\tau)u(\tau,x)\,d\tau.
\end{equation}
Assuming in addition that $u\ge 0$ on $(0,t_1)\times \Omega$ it
follows from (\ref{u0weg}), (\ref{shiftprop}), and the positivity
of $\psi$ and of $-\dot{g}_{1-\alpha,n}$ that
\begin{equation} \label{shiftprob}
\int_\Omega \Big(\psi \partial_s(g_{1-\alpha,n}\ast
\tilde{u})+\big((h_n\ast [ADu])\,\tilde{}\;|D\psi\big)\Big)\,dx \ge
\,0,\quad\mbox{a.a.}\;s\in (0,T-t_1),\,n\in \iN,
\end{equation}
for any nonnegative function $\psi\in \oH^1_2(\Omega)$. This
relation will be the starting point for all of the estimates below.

We conclude this section with a remark on the scaling properties of
equation (\ref{MProb}). Let $t_0,r>0$ and $x_0\in \iR^N$. Suppose
$u\in Z_\alpha$ is a weak solution (subsolution, supersolution) of
(\ref{MProb}) in $(0,t_0 r^{2/\alpha})\times B(x_0,r)$. Changing the
coordinates according to $s=t/r^{2/\alpha}$ and $y=(x-x_0)/r$ and
setting $v(s,y)=u(sr^{2/\alpha},x_0+yr)$, $v_0(y)=u_0(x_0+yr)$, and
$\tilde{A}(s,y)=A(sr^{2/\alpha},x_0+yr)$, the problem for $u$ is
transformed to a problem for $v$ in $(0,t_0)\times B(0,1)$, namely
there holds with $D=D_y$ (also in the weak sense)
\begin{equation} \label{MProbScal}
\partial_s^\alpha (v-v_0)-\mbox{div}\,\big(\tilde{A}(s,y)Dv\big)=\,(\le,\,\ge)\,0,\quad
s\in (0,t_0),\,y\in B(0,1).
\end{equation}

\subsection{Mean value inequalities}
For $\sigma>0$ we put $\sigma B(x,r):=B(x,\sigma r)$. Recall that
$\mu_N$ denotes the Lebesgue measure in $\iR^N$.
\begin{satz} \label{superest1}
Let $\alpha\in(0,1)$, $T>0$, and $\Omega\subset \iR^N$ be a bounded
domain. Suppose the assumptions (H1)--(H3) are satisfied. Let
$\eta>0$ and $\delta\in (0,1)$ be fixed. Then for any $t_0\in(0,T]$
and $r>0$ with $t_0-\eta r^{2/\alpha}\ge 0$, any ball
$B=B(x_0,r)\subset\Omega$, and any weak supersolution $u\ge
\varepsilon>0$ of (\ref{MProb}) in $(0,t_0)\times B$ with $u_0\ge 0$
in $B$ , there holds
\[
\esup_{U_{\sigma'}}{u^{-1}} \le \Big(\frac{C \mu_{N+1}(U_1)^{-1}
}{(\sigma-\sigma')^{\tau_0}}\Big)^{1/\gamma}
|u^{-1}|_{L_{\gamma}(U_\sigma)},\quad \delta\le \sigma'<\sigma\le
1,\; \gamma\in (0,1].
\]
Here $U_\sigma=(t_0-\sigma\eta r^{2/\alpha},t_0)\times \sigma B$,
$0<\sigma\le 1$, $C=C(\nu,\Lambda,\delta,\eta,\alpha,N)$ and
$\tau_0=\tau_0(\alpha,N)$.
\end{satz}
{\em Proof:} We may assume that $r=1$ and $x_0=0$. In fact, in the
general case we change coordinates as $t\rightarrow
t/r^{2/\alpha}$ and $x\rightarrow (x-x_0)/r$, thereby transforming
the equation to a problem of the same type on
$(0,t_0/r^{2/\alpha})\times B(0,1)$, cf. Section \ref{SSS}.

Fix $\sigma'$ and $\sigma$ such that $\delta\le \sigma'<\sigma\le
1$ and put $B_1=\sigma B$. For $\rho\in (0,1]$ we set
$V_\rho=U_{\rho\sigma}$. Given $0<\rho'<\rho\le 1$, let
$t_1=t_0-\rho\sigma\eta $ and $t_2=t_0-\rho'\sigma\eta
 $. Then $0\le t_1<t_2<t_0$. We introduce further the
shifted time ${s}=t-t_1$ and set $\tilde{f}(s)=f(s+t_1)$, $s\in
(0,t_0-t_1)$, for functions $f$ defined on $(t_1,t_0)$. Since
$u_0\ge 0$ in $B$ and $u$ is a positive weak supersolution of
(\ref{MProb}) in $(0,t_0)\times B$, we have (cf. (\ref{shiftprob}))
\begin{equation} \label{sup0}
\int_B \Big(v \partial_s(g_{1-\alpha,n}\ast \tilde{u})+\big((h_n\ast
[ADu])\,\tilde{}\;|Dv\big)\Big)\,dx \ge \,0,\quad\mbox{a.a.}\;s\in
(0,t_0-t_1),\,n\in \iN,
\end{equation}
for any nonnegative function $v\in \oH^1_2(B)$. For $s\in
(0,t_0-t_1)$ we choose the test function $v=\psi^2
\tilde{u}^{\beta}$ with $\beta<-1$ and $\psi\in C^1_0(B_1)$ so that
$0\le \psi\le 1$, $\psi=1$ in $\rho'B_1$, supp$\,\psi\subset \rho
B_1$, and $|D \psi|\le 2/[\sigma (\rho-\rho')]$. By the fundamental
identity (\ref{fundidentity}) applied to $k=g_{1-\alpha,n}$ and the
convex function $H(y)=-(1+\beta)^{-1}y^{1+\beta}$, $y>0$, there
holds for a.a. $(s,x)\in (0,t_0-t_1)\times B$
\begin{align}
 -\tilde{u}^{\beta}\partial_{s}(g_{1-\alpha,n}\ast \tilde{u}) & \ge -\,\frac{1}{1+\beta}\,\partial_{s}
(g_{1-\alpha,n}\ast\tilde{u}^{1+\beta})+\Big(\frac{\tilde{u}^{1+\beta}}{1+\beta}\,-\tilde{u}^{1+\beta}\Big)g_{1-\alpha,n}\nonumber\\
 & =
-\,\frac{1}{1+\beta}\,\partial_{s}
(g_{1-\alpha,n}\ast\tilde{u}^{1+\beta})-\,\frac{\beta}{1+\beta}\,\tilde{u}^{1+\beta}
g_{1-\alpha,n}. \label{sup1}
\end{align}
We further have
\[ Dv=2\psi D\psi \,\tilde{u}^{\beta}+\beta\psi^2 \tilde{u}^{\beta-1}D \tilde{u}.\]
Using this and (\ref{sup1}) it follows from (\ref{sup0}) that for
a.a. $s\in (0,t_0-t_1)$
\begin{align}
-\,\frac{1}{1+\beta}\,& \int_{B_1}\psi^2\partial_{s}
(g_{1-\alpha,n}\ast\tilde{u}^{1+\beta})\,dx+|\beta|\int_{B_1}\big((h_n\ast
[ADu])\,\tilde{}\;|\psi^2 \tilde{u}^{\beta-1}D
\tilde{u}\big)\,dx \nonumber\\
\le  & \,2\int_{B_1}\big((h_n\ast [ADu])\,\tilde{}\;|\psi D\psi
\,\tilde{u}^{\beta}\big)\,dx+\,\frac{\beta}{1+\beta}\,\int_{B_1}\psi^2\tilde{u}^{1+\beta}
g_{1-\alpha,n}\,dx. \label{sup2}
\end{align}
Next, choose $\varphi\in C^1([0,t_0-t_1])$ such that $0\le
\varphi\le 1$, $\varphi=0$ in $[0,(t_2-t_1)/2]$, $\varphi=1$ in
$[t_2-t_1,t_0-t_1]$, and $0\le \dot{\varphi}\le 4/(t_2-t_1)$.
Multiplying (\ref{sup2}) by $-(1+\beta)>0$ and by $\varphi(s)$,
and convolving the resulting inequality with $g_\alpha$ yields
\begin{align}
\int_{B_1} & g_\alpha\ast
\big(\varphi\partial_{s}(g_{1-\alpha,n}\ast
[\psi^2\tilde{u}^{1+\beta}])\big)\,dx+\beta(1+\beta)\,g_\alpha\ast\int_{B_1}\big((h_n\ast
[ADu])\,\tilde{}\;|\psi^2 \tilde{u}^{\beta-1}D
\tilde{u}\big)\varphi\,dx \nonumber\\
\le \, & \,2|1+\beta|\,g_\alpha\ast\int_{B_1}\big((h_n\ast
[ADu])\,\tilde{}\;|\psi D\psi
\,\tilde{u}^{\beta}\big)\varphi\,dx+|\beta|\,g_\alpha\ast\int_{B_1}\psi^2\tilde{u}^{1+\beta}
g_{1-\alpha,n}\varphi\,dx, \label{sup3}
\end{align}
for a.a. $s\in (0,t_0-t_1)$. By Lemma \ref{comm},
\begin{align}
\int_{B_1} g_\alpha\ast &
\big(\varphi\partial_{s}(g_{1-\alpha,n}\ast
[\psi^2\tilde{u}^{1+\beta}])\big)\,dx \ge \int_{B_1} \varphi
g_\alpha\ast \big(\partial_{s}(g_{1-\alpha,n}\ast
[\psi^2\tilde{u}^{1+\beta}])\big)\,dx\nonumber\\
& -\int_0^s g_\alpha(s-\sigma)\dot{\varphi}(\sigma)
\big(g_{1-\alpha,n}\ast
\int_{B_1}\psi^2\tilde{u}^{1+\beta}\,dx\big)(\sigma)\,d\sigma.
\label{sup4}
\end{align}
Furthermore, by virtue of
\[
g_{1-\alpha,n}\ast [\psi^2\tilde{u}^{1+\beta}]\in
{}_0H^1_1([0,t_0-t_1];L_1(B_1))
\]
and $g_{1-\alpha,n}=g_{1-\alpha}\ast h_n$ as well as $g_\alpha\ast
g_{1-\alpha}=1$ we have
\begin{equation} \label{sup5}
g_\alpha\ast \partial_{s}(g_{1-\alpha,n}\ast
[\psi^2\tilde{u}^{1+\beta}])=\partial_s(g_\alpha\ast
g_{1-\alpha,n}\ast [\psi^2\tilde{u}^{1+\beta}])=h_n\ast
(\psi^2\tilde{u}^{1+\beta}).
\end{equation}
Combining (\ref{sup3}), (\ref{sup4}), and (\ref{sup5}), sending
$n\to \infty$, and selecting an appropriate subsequence, if
necessary, we thus obtain
\begin{align}
& \int_{B_1}\varphi\psi^2\tilde{u}^{1+\beta}\,dx+
\beta(1+\beta)\,g_\alpha\ast\int_{B_1}\big(\tilde{A}D\tilde{u}|\psi^2
\tilde{u}^{\beta-1}D \tilde{u}\big)\varphi\,dx\nonumber\\
\le \, &
\,2|1+\beta|\,g_\alpha\ast\int_{B_1}\big(\tilde{A}D\tilde{u}|\psi
D\psi
\,\tilde{u}^{\beta}\big)\varphi\,dx+|\beta|\,g_\alpha\ast\int_{B_1}\psi^2\tilde{u}^{1+\beta}
g_{1-\alpha}\varphi\,dx \nonumber\\
& +\int_0^s g_\alpha(s-\sigma)\dot{\varphi}(\sigma)
\big(g_{1-\alpha}\ast
\int_{B_1}\psi^2\tilde{u}^{1+\beta}\,dx\big)(\sigma)\,d\sigma,
\;\;\mbox{a.a.}\;s\in(0,t_0-t_1).
\label{sup6}
\end{align}
Put $w=\tilde{u}^{\frac{\beta+1}{2}}$. Then $Dw=\frac{\beta+1}{2}
\tilde{u}^{\frac{\beta-1}{2}} D\tilde{u}$. By assumption (H2), we
have
\begin{align}
\beta(1+\beta)\,g_\alpha\ast\int_{B_1}\big(\tilde{A}D\tilde{u}|\psi^2
\tilde{u}^{\beta-1}D \tilde{u}\big)\varphi\,dx & \,\ge \nu
\beta(1+\beta)\,g_\alpha\ast\int_{B_1} \varphi
\psi^2\tilde{u}^{\beta-1}|D\tilde{u}|^2\,dx \nonumber\\
& \, = \,\frac{4\nu \beta}{1+\beta}\,g_\alpha\ast\int_{B_1}\varphi
\psi^2|Dw|^2\,dx. \label{sup7}
\end{align}
Using (H1) and Young's inequality we may estimate
\begin{align}
2\big|\big(\tilde{A}D\tilde{u}|\psi D\psi
\,\tilde{u}^{\beta}\big)\varphi\big| & \le 2\Lambda\psi|D\psi|\,|D
\tilde{u}|\tilde{u}^\beta \varphi=2\Lambda\psi|D\psi|\,|D
\tilde{u}|\tilde{u}^{\frac{\beta-1}{2}}\tilde{u}^{\frac{\beta+1}{2}}\varphi\nonumber\\
& \le \,\frac{\nu |\beta|}{2}\, \psi^2\varphi |D \tilde{u}|^2
\tilde{u}^{\beta-1}+\,\frac{2}{\nu |\beta|}\,\Lambda^2
|D\psi|^2\varphi \tilde{u}^{\beta+1}\nonumber\\
& = \,\frac{2\nu |\beta|}{(1+\beta)^2}\,
\psi^2\varphi|Dw|^2+\,\frac{2}{\nu |\beta|}\,\Lambda^2
|D\psi|^2\varphi w^2. \label{sup8}
\end{align}
From (\ref{sup6}), (\ref{sup7}), and (\ref{sup8}) we conclude that
\begin{equation} \label{sup9}
\int_{B_1}\varphi\psi^2w^2\,dx+\,\frac{2\nu
|\beta|}{|1+\beta|}\,g_\alpha\ast\int_{B_1}\varphi \psi^2|Dw|^2\,dx
\le g_\alpha\ast F,\quad\mbox{a.a.}\;s\in(0,t_0-t_1),
\end{equation}
where
\begin{align*}
F(s) =\, & \,\frac{2\Lambda^2|1+\beta|}{\nu |\beta|}\, \int_{B_1}
|D\psi|^2\varphi w^2\,dx
+|\beta|\varphi(s)g_{1-\alpha}(s)\int_{B_1}\psi^2 w^2
\,dx \\
& \,+\dot{\varphi}(s) \big(g_{1-\alpha}\ast \int_{B_1}\psi^2
w^2\,dx\big)(s)\ge 0,\quad\mbox{a.a.}\;s\in(0,t_0-t_1).
\end{align*}
We may drop the second term in (\ref{sup9}), which is nonnegative.
By Young's inequality for convolutions and the properties of
$\varphi$ we then infer that for all $p\in(1,1/(1-\alpha))$
\begin{equation} \label{sup10}
\Big(\int_{t_2-t_1}^{t_0-t_1} (\int_{B_1}
[\psi(x)w(s,x)]^2\,dx)^p\,ds\Big)^{1/p} \,\le
|g_\alpha|_{L_p([0,t_0-t_1])} \int_0^{t_0-t_1} \!\!\!\!F(s)\,ds,
\end{equation}
where
\begin{equation} \label{sup11}
|g_\alpha|_{L_p([0,t_0-t_1])}=
\,\frac{(t_0-t_1)^{\alpha-1+1/p}}{\Gamma(\alpha)[(\alpha-1)p+1]^{1/p}}\,\le
\,\frac{\eta^{\alpha-1+1/p}}{\Gamma(\alpha)[(\alpha-1)p+1]^{1/p}}\,
=:C_1(\alpha,p,\eta).
\end{equation}
We choose any of these $p$ and fix it.

Returning to (\ref{sup9}), we may also drop the first term, convolve
the resulting inequality with $g_{1-\alpha}$ and evaluate at
$s=t_0-t_1$, thereby obtaining
\begin{equation} \label{sup12}
\int_{t_2-t_1}^{t_0-t_1}\int_{B_1}\psi^2|Dw|^2\,dx\,ds \le
\,\frac{|1+\beta|}{2\nu |\beta|}\,\int_0^{t_0-t_1} \!\!\!\!F(s)\,ds.
\end{equation}
Using
\[
\int_{t_2-t_1}^{t_0-t_1}\int_{B_1} |D(\psi w)|^2\,dx\,ds\le
2\int_{t_2-t_1}^{t_0-t_1}\int_{B_1}
\big(\psi^2|Dw|^2+|D\psi|^2w^2\big)\,dx\,ds
\]
we infer from (\ref{sup10})--(\ref{sup12}) that
\begin{align}
|\psi w|^2_{V_p([t_2-t_1,t_0-t_1]\times B_1)}\le &\;
2\Big(C_1(\alpha,p,\eta)+\,\frac{|1+\beta|}{\nu|\beta|}\Big)\int_0^{t_0-t_1}
\!\!\!\!F(s)\,ds\nonumber\\
&+4\int_{0}^{t_0-t_1}\int_{B_1} |D\psi|^2w^2\,dx\,ds.
\label{vpest}
\end{align}

We will next estimate the right-hand side of (\ref{vpest}). By the
assumptions on $\psi$ and $\varphi$, and since $|\beta|>1$, we
have
\[ \int_0^{t_0-t_1}\!\!\!\int_{B_1} |D\psi|^2
w^2\,dx\,ds\le
\,\frac{4}{\sigma^2(\rho-\rho')^2}\,\int_0^{t_0-t_1}
\!\!\!\int_{\rho B_1}w^2\,dx\,d{s}
\]
and
\begin{align*}
F(s)\le &\;\Big(\,\frac{8\Lambda^2|1+\beta|}{\nu
\sigma^2(\rho-\rho')^2}\,+|\beta|g_{1-\alpha}((t_2-t_1)/2)\Big)\int_{\rho
B_1}w^2\,dx \nonumber\\ &
+\,\frac{4}{t_2-t_1}\,\big(g_{1-\alpha}\ast \int_{\rho B_1}
w^2\,dx\big)(s),\quad \mbox{a.a.}\;s\in (0,t_0-t_1).
\end{align*}
Recall that $\sigma\ge \delta>0$. So we have
\begin{align*}
\int_0^{t_0-t_1} \!\!\!\!F(s)\,ds \;\le &\;
\Big(\,\frac{8\Lambda^2|1+\beta|}{\nu
\sigma^2(\rho-\rho')^2}\,+\,\frac{2^\alpha|\beta|}
{\Gamma(1-\alpha)(\rho-\rho')^\alpha(\sigma\eta)^\alpha}
\Big)\int_0^{t_0-t_1} \!\!\!\int_{\rho B_1}w^2\,dx\,ds\\
&
\;+\,\frac{4}{(\rho-\rho')\sigma\eta}\,\int_0^{t_0-t_1}g_{2-\alpha}(t_0-t_1-\tau)
\int_{\rho B_1}w(\tau,x)^2\,dx\,d\tau \\
\le & \;
C(\nu,\Lambda,\delta,\eta,\alpha)\,\frac{1+|1+\beta|}{(\rho-\rho')^2}
\int_0^{t_0-t_1} \!\!\!\int_{\rho B_1}w^2\,dx\,ds.
\end{align*}
Combining these estimates and (\ref{vpest}) yields
\[
|\psi w|_{V_p([t_2-t_1,t_0-t_1]\times B_1)}\le
C(\nu,\Lambda,\delta,\eta,\alpha,p)\,\frac{1+|1+\beta|}
{\rho-\rho'}\, |w|_{L_{2}([0,t_0-t_1]\times\rho B_1)}.
\]
We apply next the interpolation inequality (\ref{Vembedding}) to
the function $\psi w$ and make use of $\psi=1$ in $\rho'B_1$ to
deduce that
\begin{equation} \label{sup13}
|w|_{L_{2\kappa}([t_2-t_1,t_0-t_1]\times \rho'B_1)}\le
C(\nu,\Lambda,\delta,\eta,\alpha,p,N)\,\frac{1+|1+\beta|}{\rho-\rho'}
\,|w|_{L_{2}([0,t_0-t_1]\times \rho B_1)},
\end{equation}
where the number $\kappa>1$ is given in (\ref{kappa}). Since
$w=\tilde{u}^{\frac{\beta+1}{2}}$ and by transforming back to the
time $t$, we see that (\ref{sup13}) is equivalent to
\[
(\int_{V_{\rho'}}
u^{-|1+\beta|\kappa}\,d\mu_{N+1})^{\frac{1}{2\kappa}}\le
\frac{\tilde{C}(1+|1+\beta|)}{\rho-\rho'}\,(\int_{V_{\rho}}
u^{-|1+\beta|}\,d\mu_{N+1})^{\frac{1}{2}}
\]
with $\tilde{C}=\tilde{C}(\nu,\Lambda,\delta,\eta,\alpha,p,N)$.
Hence, with $\gamma=|1+\beta|$,
\[
|u^{-1}|_{L_{\gamma\kappa}(V_{\rho'})}\le
\Big(\frac{\tilde{C}^2(1+\gamma)^2}{(\rho-\rho')^2}\Big)^{1/\gamma}
|u^{-1}|_{L_{\gamma}(V_{\rho})},\quad 0<\rho'<\rho\le
1,\;\gamma>0.
\]
Employing the first Moser iteration, Lemma \ref{moserit1} (with
$\bar{p}=1$), it follows that there exist constants
$M_0=M_0(\nu,\Lambda,\delta,\eta,\alpha,p,N)$ and
$\tau_0=\tau_0(\kappa)$ such that
\[
\esup_{V_{\theta}}{u^{-1}} \le
\Big(\frac{M_0}{(1-\theta)^{\tau_0}}\Big)^{1/\gamma}
|u^{-1}|_{L_{\gamma}(V_1)}\quad \mbox{for
all}\;\;\theta\in(0,1),\;\gamma\in (0,1] .\] Thus if we take
$\theta=\sigma'/\sigma$ and notice that
\[ \frac{1}{1-\theta}\,=\,\frac{\sigma}{\sigma-\sigma'}\,\le
\frac{1}{\sigma-\sigma'},\] we obtain
\[
\esup_{U_{\sigma'}}{u^{-1}} \le
\Big(\frac{M_0}{(\sigma-\sigma')^{\tau_0}}\Big)^{1/\gamma}
|u^{-1}|_{L_{\gamma}(U_\sigma)},\quad \gamma\in (0,1].
\]
Hence the proof is complete. $\square$

${}$

We put
\[
\tilde{\kappa}:=\kappa_{1/(1-\alpha)}=\,\frac{2+N\alpha}{2+N\alpha-2\alpha}.
\]
\begin{satz} \label{superest2}
Let $\alpha\in(0,1)$, $T>0$, and $\Omega\subset \iR^N$ be a bounded
domain. Suppose the assumptions (H1)--(H3) are satisfied. Let
$\eta>0$ and $\delta\in (0,1)$ be fixed. Then for any $t_0\in [0,T)$
and $r>0$ with $t_0+\eta r^{2/\alpha}\le T$, any ball
$B=B(x_0,r)\subset\Omega$, any $p_0\in(0,\tilde{\kappa})$, and any
nonnegative weak supersolution $u$ of (\ref{MProb}) in $(0,t_0+\eta
r^{2/\alpha})\times B$ with $u_0\ge 0$ in $B$, there holds
\[
|u|_{L_{p_0}(U_{\sigma'}')}\le \Big(\frac{C\mu_{N+1}(U'_1)^{-1}
}{(\sigma-\sigma')^{\tau_0}}\Big)^{1/\gamma-1/p_0}
|u|_{L_{\gamma}(U'_\sigma)},\quad \delta\le \sigma'<\sigma\le 1,\;
0<\gamma\le p_0/\tilde{\kappa}.
\]
Here $U'_\sigma=(t_0,t_0+\sigma\eta r^{2/\alpha})\times \sigma B$,
$C=C(\nu,\Lambda,\delta,\eta,\alpha,N,p_0)$, and
$\tau_0=\tau_0(\alpha,N)$.
\end{satz}
{\em Proof:} We proceed similarly as in the previous proof.
Without restriction of generality we may assume that $p_0>1$ and
$r=1$. By replacing $u$ with $u+\varepsilon$ and $u_0$ with
$u_0+\varepsilon$ and eventually letting $\varepsilon \to 0+$ we
may further assume that $u$ is bounded away from zero.

Fix $\sigma'$, $\sigma$ such that $\delta\le \sigma'<\sigma\le 1$
and put $B_1=\sigma B$. For $\rho\in (0,1]$ we set
$V'_\rho=U'_{\rho\sigma}$. Given $0<\rho'<\rho\le 1$, let
$t_1=t_0+\rho'\sigma\eta$ and $t_2=t_0+\rho\sigma\eta$, so $0\le
t_0<t_1<t_2$. We shift the time by means of ${s}=t-t_0$ and set
$\tilde{f}(s)=f(s+t_0)$, $s\in (0,t_2-t_0)$, for functions $f$
defined on $(t_0,t_2)$.

We then repeat the first steps of the preceding proof, the only
difference being that now we take $\beta\in (-1,0)$. Note that, as
a consequence of this, (\ref{sup1}) simplifies to
\[
-\tilde{u}^{\beta}\partial_{s}(g_{1-\alpha,n}\ast \tilde{u}) \ge
-\,\frac{1}{1+\beta}\,\partial_{s}
(g_{1-\alpha,n}\ast\tilde{u}^{1+\beta}),\quad
\mbox{a.a.}\;(s,x)\in (0,t_2-t_0)\times B,
\]
hence we obtain with $\psi\in C_0^1(B_1)$ as above
\begin{align}
-\,\frac{1}{1+\beta}\,& \int_{B_1}\psi^2\partial_{s}
(g_{1-\alpha,n}\ast\tilde{u}^{1+\beta})\,dx+|\beta|\int_{B_1}\big((h_n\ast
[ADu])\,\tilde{}\;|\psi^2 \tilde{u}^{\beta-1}D
\tilde{u}\big)\,dx \nonumber\\
\le  & \,2\int_{B_1}\big((h_n\ast [ADu])\,\tilde{}\;|\psi D\psi
\,\tilde{u}^{\beta}\big)\,dx,\quad\quad
\mbox{a.a.}\;s\in(0,t_2-t_0). \label{L1}
\end{align}

Next, choose $\varphi\in C^1([0,t_2-t_0])$ such that $0\le
\varphi\le 1$, $\varphi=1$ in $[0,t_1-t_0]$, $\varphi=0$ in
$[t_1-t_0+(t_2-t_1)/2,t_2-t_0]$, and $0\le -\dot{\varphi}\le
4/(t_2-t_1)$. Multiplying (\ref{L1}) by $1+\beta>0$ and by
$\varphi(s)$, and applying Lemma \ref{comm2} to the first term
gives
\begin{align}
-\int_{B_1}  &
\partial_{s}(g_{1-\alpha,n}\ast
[\varphi\psi^2\tilde{u}^{1+\beta}]\big)\,dx+|\beta|(1+\beta)\,
\int_{B_1}\big(\tilde{A}D\tilde{u}|\psi^2 \tilde{u}^{\beta-1}D
\tilde{u}\big)\varphi\,dx \nonumber\\
\le & \,\int_0^s
\dot{g}_{1-\alpha,n}(s-\sigma)\big(\varphi(s)-\varphi(\sigma)\big)
\big(\int_{B_1}\psi^2\tilde{u}^{1+\beta}\,dx\big)(\sigma)\,d\sigma\nonumber\\
& \;\,+2(1+\beta)\,\int_{B_1}\big(\tilde{A}D\tilde{u}|\psi D\psi
\,\tilde{u}^{\beta}\big)\varphi\,dx+\mathcal{R}_n(s) ,\quad
\mbox{a.a.}\;s\in(0,t_2-t_0), \label{L2}
\end{align}
where
\begin{align*}
\mathcal{R}_n(s)= &\,\,-|\beta|(1+\beta)\, \int_{B_1}\big((h_n\ast
[ADu])\,\tilde{}\;-\tilde{A}D\tilde{u}|\psi^2 \tilde{u}^{\beta-1}D
\tilde{u}\big)\varphi\,dx\\
&\,+2(1+\beta)\,\int_{B_1}\big((h_n\ast
[ADu])\,\tilde{}\;-\tilde{A}D\tilde{u}|\psi D\psi
\,\tilde{u}^{\beta}\big)\varphi\,dx,\quad
\mbox{a.a.}\;s\in(0,t_2-t_0).
\end{align*}
We set again $w=\tilde{u}^{\frac{\beta+1}{2}}$ and estimate
exactly as in the preceding proof, using (H1), (H3) and
(\ref{sup8}), to the result
\begin{align}
-\int_{B_1}  &
\partial_{s}(g_{1-\alpha,n}\ast
[\varphi\psi^2w^2]\big)\,dx+\,\frac{2\nu
|\beta|}{1+\beta}\,\int_{B_1}\varphi
\psi^2|Dw|^2\,dx \nonumber\\
\le & \,\int_0^s
\dot{g}_{1-\alpha,n}(s-\sigma)\big(\varphi(s)-\varphi(\sigma)\big)
\big(\int_{B_1}\psi^2w^2\,dx\big)(\sigma)\,d\sigma\nonumber\\
& \;\,+\,\frac{2\Lambda^2(1+\beta)}{\nu |\beta|}\, \int_{B_1}
|D\psi|^2\varphi w^2\,dx+\mathcal{R}_n(s) ,\quad
\mbox{a.a.}\;s\in(0,t_2-t_0). \label{L3}
\end{align}
Recall that $g_{1-\alpha,n}=g_{1-\alpha}\ast h_n$. Putting
\[
W(s)=\int_{B_1}\varphi(s)\psi(x)^2w(s,x)^2\,dx
\]
and denoting the right-hand side of (\ref{L3}) by $F_n(s)$, it
follows from (\ref{L3}) that
\[
G_n(s):=\partial_s^\alpha (h_n\ast W)(s)+F_n(s)\ge 0,\quad\quad
\mbox{a.a.}\;s\in(0,t_2-t_0).
\]
By (\ref{sup5}) and positivity of $h_n$, we have
\[
0\le h_n\ast W =g_\alpha\ast \partial_s^\alpha (h_n\ast W)\le
g_\alpha\ast G_n+g_\alpha\ast [-F_n(s)]_+
\]
a.e. in $(0,t_2-t_0)$, where $[y]_+$ stands for the positive part
of $y\in \iR$. For any $p\in (1,1/(1-\alpha))$ and any
$t_*\in[t_2-t_0-(t_2-t_1)/4,t_2-t_0]$ we thus obtain by Young's
inequality
\begin{equation} \label{L4}
|h_n\ast W|_{L_p([0,t_*])}\le
|g_\alpha|_{L_p([0,t_*])}\big(|G_n|_{L_1([0,t_*])}+
|[-F_n]_+|_{L_1([0,t_*])}\big).
\end{equation}
Since $t_*\le t_2-t_0\le \eta$, we have
$|g_\alpha|_{L_p([0,t_*])}\le C_1(\alpha,p,\eta)$ with the same
constant as in (\ref{sup11}). By positivity of $G_n$,
\[
|G_n|_{L_1([0,t_*])}=(g_{1-\alpha,n}\ast
W)(t_*)+\int_0^{t_*}\!\!\!F_n(s)\,ds.
\]
Observe that $\mathcal{R}_n\rightarrow 0$ in $L_1([0,t_2-t_0])$ as
$n\to \infty$. Hence $|[-F_n]_+|_{L_1([0,t_*])}\to 0$ as
$n\to\infty$. Further,
\begin{align*}
\int_0^{t_*}\!\!&\!\int_0^s
\dot{g}_{1-\alpha,n}(s-\sigma)\big(\varphi(s)-\varphi(\sigma)\big)
\big(\int_{B_1}\psi^2w^2\,dx\big)(\sigma)\,d\sigma\,ds\\
=&\,\int_0^{t_*}
{g}_{1-\alpha,n}(t_*-\sigma)\big(\varphi(t_*)-\varphi(\sigma)\big)
\big(\int_{B_1}\psi^2w^2\,dx\big)(\sigma)\,d\sigma\\
&\,-\int_0^{t_*}\!\!\!\dot{\varphi}(s)\int_0^s
{g}_{1-\alpha,n}(s-\sigma)
\big(\int_{B_1}\psi^2w^2\,dx\big)(\sigma)\,d\sigma\,ds\\
\le&\,-\int_0^{t_*}\!\!\!\dot{\varphi}(s)\int_0^s
{g}_{1-\alpha,n}(s-\sigma)
\big(\int_{B_1}\psi^2w^2\,dx\big)(\sigma)\,d\sigma\,ds,
\end{align*}
since $\varphi$ is nonincreasing. We also know that
$g_{1-\alpha,n}\ast W\to g_{1-\alpha}\ast W$ in
$L_1([0,t_2-t_0])$. Hence we can fix some
$t_*\in[t_2-t_0-(t_2-t_1)/4,t_2-t_0]$ such that for some
subsequence $(g_{1-\alpha,n_k}\ast W)(t_*)\to (g_{1-\alpha}\ast
W)(t_*)$ as $k\to \infty$. Sending $k\to \infty$ it follows then
from (\ref{L4}), the preceding estimates, and from $\varphi=1$ in
$[0,t_1-t_0]$ that
\begin{equation} \label{L5}
\Big(\int_{0}^{t_1-t_0} (\int_{B_1}
[\psi(x)w(s,x)]^2\,dx)^p\,ds\Big)^{1/p}\le
C_1(\alpha,p,\eta)\Big((g_{1-\alpha}\ast
W)(t_*)+|F|_{L_1([0,t_2-t_0])}\Big),
\end{equation}
with
\[
F(s)=\,\frac{2\Lambda^2(1+\beta)}{\nu |\beta|}\, \int_{B_1}
|D\psi|^2\varphi w^2\,dx-\dot{\varphi}(s)\big(g_{1-\alpha}\ast
\int_{B_1}\psi^2w^2\,dx\big)(s).
\]

On the other hand, we can integrate (\ref{L3}) over $(0,t_*)$ and
take the limit as $k\to \infty$ for the same subsequence as
before, thereby getting
\begin{equation} \label{L6}
\int_{0}^{t_1-t_0}\!\!\!\int_{B_1}
\psi^2|Dw|^2\,dx\,ds\le\,\frac{1+\beta}{2\nu
|\beta|}\,\Big((g_{1-\alpha}\ast
W)(t_*)+|F|_{L_1([0,t_2-t_0])}\Big).
\end{equation}
Arguing as above (cf. the lines before (\ref{vpest})), we conclude
from (\ref{L5}) and (\ref{L6}) that
\begin{align}
&|\psi w|^2_{V_p([0,t_1-t_0]\times B_1)}\le \;
4\int_{0}^{t_2-t_0}\int_{B_1} |D\psi|^2w^2\,dx\,ds\nonumber\\
&+2\Big(C_1(\alpha,p,\eta)+\,\frac{1+\beta}{\nu|\beta|}\Big)
\Big((g_{1-\alpha}\ast W)(t_*)+|F|_{L_1([0,t_2-t_0])}\Big).
\label{L7}
\end{align}

Since $\varphi=0$ in $[t_1-t_0+(t_2-t_1)/2,t_2-t_0]$ and
$t_*\in[t_2-t_0-(t_2-t_1)/4,t_2-t_0]$, we have
\begin{align*}
(g_{1-\alpha}\ast W)(t_*)&\,\le
g_{1-\alpha}\big((t_2-t_1)/4\big)\int_0^{t_2-t_0} \!\!\!\int_{
B_1}\psi^2w^2\,dx\,d{s}\\
&\,=\,\frac{4^\alpha}{\Gamma(1-\alpha)(\rho-\rho')^\alpha(\sigma\eta)^\alpha}\,
\int_0^{t_2-t_0} \!\!\!\int_{\rho B_1}w^2\,dx\,d{s}.
\end{align*}
Further,
\[
 \int_{0}^{t_2-t_0}\int_{B_1} |D\psi|^2w^2\,dx\,ds\le\,\frac{4}{\sigma^2(\rho-\rho')^2}
 \,\int_0^{t_2-t_0}
\!\!\!\int_{\rho B_1}w^2\,dx\,d{s}.
\]
The term $|F|_{L_1([0,t_2-t_0])}$ is estimated similarly as in the
proof of Theorem \ref{superest1} (cf. the lines that follow
(\ref{vpest})). We obtain
\begin{align*}
|F|_{L_1([0,t_2-t_0])}\le
C(\nu,\Lambda,\delta,\eta,\alpha)\,\frac{1+(1+\beta)}{|\beta|(\rho-\rho')^2}
\int_0^{t_2-t_0} \!\!\!\int_{\rho B_1}w^2\,dx\,ds.
\end{align*}
Notice the additional factor $|\beta|$ in the denominator.
Combining these estimates we deduce from (\ref{L7}) that
\[
|\psi w|_{V_p([0,t_1-t_0]\times B_1)}\le
C(\nu,\Lambda,\delta,\eta,\alpha,p)\,\frac{1+(1+\beta)}
{|\beta|(\rho-\rho')}\, |w|_{L_{2}([0,t_2-t_0]\times\rho B_1)}.
\]
By the interpolation inequality (\ref{Vembedding}) and since
$\psi=1$ in $\rho'B_1$, this implies for all $\beta\in(-1,0)$
\begin{equation} \label{supL8}
|w|_{L_{2\kappa}([0,t_1-t_0]\times \rho'B_1)}\le
C(\nu,\Lambda,\delta,\eta,\alpha,p,N)\,\frac{1+|1+\beta|}{|\beta|(\rho-\rho')}
\,|w|_{L_{2}([0,t_2-t_0]\times \rho B_1)},
\end{equation}
where
\[
\kappa=\kappa_p=\,\frac{2p+N(p-1)}{2+N(p-1)}\,\in
(1,\tilde{\kappa}).
\]
We now fix $1<p<1/(1-\alpha)$ such that
$\kappa_p=(p_0+\tilde{\kappa})/2$. This is possible because
$\kappa_p\nearrow \tilde{\kappa}$ as $p \nearrow 1/(1-\alpha)$.

Next, we set $\gamma=1+\beta\in (0,1)$ and transform back to $u$
to get
\begin{equation} \label{L9}
|u|_{L_{\gamma\kappa}(V'_{\rho'},d\mu)}\le
\Big(\frac{\tilde{C}}{(\rho-\rho')^{2}}\Big)^{1/\gamma}
|u|_{L_{\gamma}(V'_{\rho},d\mu)},\quad 0<\rho'<\rho\le
1,\;0<\gamma\le p_0/\kappa.
\end{equation}
Here, $\mu=(\eta \omega_N)^{-1}\mu_{N+1}$, $\omega_N$ the volume
of the unit ball in $\iR^N$, and
$\tilde{C}=\tilde{C}(\nu,\Lambda,\delta,\eta,\alpha,N,p_0)$ is
independent of $\gamma\in(0,p_0/\kappa]$, since $|\beta|$ is
bounded away from zero. Note that $\mu(V'_1)\le 1$.

Finally, we employ the second Moser iteration scheme, Lemma
\ref{moserit2}, to conclude from (\ref{L9}) that there are
constants $M_0=M_0(\nu,\Lambda,\delta,\eta,\alpha,N,p_0)$ and
$\tau_0=\tau_0(\kappa)$ such that
\begin{equation} \label{L10}
|u|_{L_{p_0}(V'_{\theta},d\mu)}\le
\Big(\frac{M_0}{(1-\theta)^{\tau_0}}\Big)^{1/\gamma-1/p_0}
|u|_{L_{\gamma}(V'_{1},d\mu)},\quad 0<\theta< 1,\;0<\gamma\le
p_0/\kappa.
\end{equation}
If we take $\theta=\sigma'/\sigma$ and translate (\ref{L10}) back
to the measure $\mu_{N+1}$, we obtain
\begin{equation} \label{L11}
|u|_{L_{p_0}(U_{\sigma'}')}\le \Big(\frac{M_0 (\eta\omega_N)^{-1}
}{(\sigma-\sigma')^{\tau_0}}\Big)^{1/\gamma-1/p_0}
|u|_{L_{\gamma}(U'_\sigma)},\quad 0<\gamma\le p_0/\kappa.
\end{equation}
Since $\kappa<\tilde{\kappa}$, (\ref{L11}) holds in particular for
all $\gamma\in (0,p_0/\tilde{\kappa}]$. This finishes the proof.
$\square$
\subsection{Logarithmic estimates}
\begin{satz} \label{logest}
Let $\alpha\in(0,1)$, $T>0$, and $\Omega\subset \iR^N$ be a bounded
domain. Suppose the assumptions (H1)--(H3) are satisfied. Let
$\tau>0$ and $\delta,\,\eta\in(0,1)$ be fixed. Then for any $t_0\ge
0$ and $r>0$ with $t_0+\tau r^{2/\alpha}\le T$, any ball
$B=B(x_0,r)\subset\Omega$, and any weak supersolution $u\ge
\varepsilon>0$ of (\ref{MProb}) in $(0,t_0+\tau r^{2/\alpha})\times
B$ with $u_0\ge 0$ in $B$, there is a constant $c=c(u)$ such that
\begin{equation} \label{logestleft}
\mu_{N+1}\big(\{(t,x)\in K_-: \log u(t,x)>c+\lambda\}\big)\le C
r^{2/\alpha} \mu_N(B)\lambda^{-1},\quad \lambda>0,
\end{equation}
and
\begin{equation} \label{logestright}
\mu_{N+1}\big(\{(t,x)\in K_+: \log u(t,x)<c-\lambda\}\big)\le C
r^{2/\alpha} \mu_N(B)\lambda^{-1},\quad \lambda>0,
\end{equation}
where $K_-=(t_0,t_0+\eta \tau r^{2/\alpha})\times \delta B$ and
$K_+=(t_0+\eta \tau r^{2/\alpha},t_0+\tau r^{2/\alpha})\times
\delta B$. Here the constant $C$ depends only on $\delta, \eta,
\tau, N, \alpha, \nu$, and $\Lambda$.
\end{satz}
{\em Proof:} Since $u_0\ge 0$ in $B$ and $u$ is a positive weak
supersolution we may assume without loss of generality that $u_0=0$
and $t_0=0$. In fact, in the case $t_0>0$ we shift the time as $t\to
t-t_0$, thereby obtaining an inequality of the same type on the
time-interval $J:=[0,\tau r^{2/\alpha}]$. Observe that the property
$g_{1-\alpha}\ast u\in C([0,t_0+\tau r^{2/\alpha}];L_2(B))$ implies
$g_{1-\alpha}\ast \tilde{u}\in C(J;L_2(B))$ for the shifted function
$\tilde{u}(s,x)=u(s+t_0,x)$. So we have
\begin{equation} \label{log1}
\int_B \Big(v \partial_t(g_{1-\alpha,n}\ast u)+(h_n\ast
[ADu]|Dv)\Big)\,dx \ge \,0,\quad\mbox{a.a.}\;t\in J,\,n\in \iN,
\end{equation}
for any nonnegative test function $v\in \oH^1_2(B)$.

For $t\in J$ we choose the test function $v=\psi^2 u^{-1}$ with
$\psi\in C^1_0(B)$ such that supp$\,\psi\subset B$, $\psi=1$ in
$\delta B$, $0\le \psi \le 1$, $|D \psi|\le 2/[(1-\delta)r]$ and
the domains $\{x\in B:\psi(x)^2 \ge b\}$ are convex for all $b\le
1$. We have
\[ Dv=2\psi D\psi \,u^{-1}-\psi^2 u^{-2}Du,\]
so that by substitution into (\ref{log1}) we obtain for a.a. $t\in
J$
\begin{align} \label{log1a}
-\int_B \psi^2 u^{-1}\partial_t & (g_{1-\alpha,n}\ast
u)\,dx+\int_B\big(ADu|u^{-2}Du\big)\psi^2\,dx\nonumber\\ &\le
2\int_B\big(ADu|u^{-1}\psi D\psi\big)\,dx+\mathcal{R}_n(t),
\end{align}
where
\[
\mathcal{R}_n(t)=\int_B\big(h_n\ast [ADu]-ADu|Dv\big)\,dx.
\]
By (H1) and Young's inequality,
\[
\big|2\big(ADu|u^{-1}\psi D\psi\big)\big|\le 2\Lambda \psi
|D\psi|\,|Du| u^{-1}\le \frac{\nu}{2}\,\psi^2|Du|^2
u^{-2}+\frac{2}{\nu}\,\Lambda^2|D\psi|^2.
\]
Using this, (H2) and $|D \psi|\le 2/[(1-\delta)r]$, we infer from
(\ref{log1a}) that for a.a. $t\in J$
\begin{equation} \label{log2}
-\int_B \psi^2 u^{-1}\partial_t(g_{1-\alpha,n}\ast
u)\,dx+\frac{\nu}{2}\,\int_B |Du|^2 u^{-2} \psi^2\,dx \le \frac{8
\Lambda^2 \mu_N(B)}{\nu (1-\delta)^2 r^2}\,+\mathcal{R}_n(t).
\end{equation}
Setting $w=\log u$ we have $Dw=u^{-1} Du$. The weighted Poincar\'e
inequality of Proposition \ref{WeiPI} with weight $\psi^2$ yields
\begin{equation} \label{log3}
\int_B (w-W)^2 \psi^2 dx \le \frac{8 r^2 \mu_N(B)}{\int_B \psi^2
dx}\,\int_B |Dw|^2 \psi^2 dx,\quad \mbox{a.a.}\;t\in J,
\end{equation}
where
\[ W(t)=\,\frac{\int_B w(t,x) \psi(x)^2 dx}{\int_B \psi(x)^2
dx}\,,\quad\quad \mbox{a.a.}\;t\in J.
\]
From (\ref{log2}) and (\ref{log3}) we deduce that
\[
-\int_B \psi^2 u^{-1}\partial_t(g_{1-\alpha,n}\ast u)
\,dx+\,\frac{\nu \int_B \psi^2 dx}{16r^2 \mu_N(B)}\,\int_B (w-W)^2
\psi^2 dx \le \frac{8 \Lambda^2 \mu_N(B)}{\nu (1-\delta)^2
r^2}\,+\mathcal{R}_n(t),
\]
which in turn implies
\begin{equation} \label{log4}
\frac{-\int_B \psi^2 u^{-1}\partial_t(g_{1-\alpha,n}\ast u)
\,dx}{\int_B \psi^2 dx}+\,\frac{\nu}{16r^2 \mu_N(B)}\,\int_{\delta
B} (w-W)^2 dx \le \frac{C_1}{r^2}\,+S_n(t),
\end{equation}
for a.a. $t\in J$, with some constant
$C_1=C_1(\delta,N,\nu,\Lambda)$ and
$S_n(t)=\mathcal{R}_n(t)/\int_B \psi^2 dx$.

The fundamental identity (\ref{fundidentity}) with $H(y)=-\log y$
reads (with the spatial variable $x$ being suppressed)
\begin{align*}
-u^{-1}\partial_t(g_{1-\alpha,n} & \ast
u)=-\partial_t(g_{1-\alpha,n}\ast \log u)+(\log
u-1)g_{1-\alpha,n}(t)\nonumber\\
&+\int_0^t \Big(-\log u(t-s)+\log
u(t)+\frac{u(t-s)-u(t)}{u(t)}\Big)[-\dot{g}_{1-\alpha,n}(s)]\,ds.
\end{align*}
In terms of $w=\log u$ this means that
\begin{align} \label{log5}
-u^{-1}\partial_t(g_{1-\alpha,n}  \ast u)= & \,
-\partial_t(g_{1-\alpha,n}\ast
w)+(w-1)g_{1-\alpha,n}(t)\nonumber\\
& \,\,+\int_0^t
\Psi\big(w(t-s)-w(t)\big)[-\dot{g}_{1-\alpha,n}(s)]\,ds,
\end{align}
where $\Psi(y)=e^y-1-y$. Since $\Psi$ is convex, it follows from
Jensen's inequality that
\[
\frac{\int_B \psi^2 \Psi\big(w(t-s,x)-w(t,x)\big)\,dx}{\int_B
\psi^2 dx} \ge \Psi \Big( \frac{\int_B \psi^2
\big(w(t-s,x)-w(t,x)\big)\,dx}{\int_B \psi^2 dx}\Big).
\]
Using this and (\ref{log5}) we obtain
\begin{align}
\frac{-\int_B \psi^2 u^{-1}\partial_t(g_{1-\alpha,n}\ast u)
\,dx}{\int_B \psi^2 dx} & \ge -\partial_t(g_{1-\alpha,n}\ast
W)+(W-1)g_{1-\alpha,n}(t)\nonumber\\
&\quad +\int_0^t
\Psi\big(W(t-s)-W(t)\big)[-\dot{g}_{1-\alpha,n}(s)]\,ds \nonumber\\
& = -e^{-W} \partial_t(g_{1-\alpha,n}\ast e^W), \label{log6}
\end{align}
where the last equals sign holds again by (\ref{log5}) with $u$
replaced by $e^W$. From (\ref{log4}) and (\ref{log6}) we conclude
that
\begin{equation} \label{log7}
\frac{\nu}{16r^2 \mu_N(B)}\,\int_{\delta B} (w-W)^2 dx \le e^{-W}
\partial_t(g_{1-\alpha,n}\ast
e^W)+\,\frac{C_1}{r^2}\,+S_n(t),\quad \mbox{a.a.}\;t\in J.
\end{equation}

We choose
\begin{equation} \label{cwahl}
c(u)=\log \Big(\frac{(g_{1-\alpha}\ast e^W)(\eta\tau
r^{2/\alpha})}{g_{2-\alpha}(\eta\tau r^{2/\alpha})}\Big).
\end{equation}
This definition makes sense, since $g_{1-\alpha}\ast e^W\in C(J)$.
The latter is a consequence of $g_{1-\alpha}\ast u\in C(J;L_2(B))$
and
\[
e^{W(t)}\le \,\frac{\int_B u(t,x) \psi(x)^2 dx}{\int_B \psi(x)^2
dx}\,,\quad\quad \mbox{a.a.}\;t\in J,
\]
where we apply again Jensen's inequality.

To prove (\ref{logestleft}) and (\ref{logestright}), one of the
key ideas is to use the inequalities
\begin{align}
\mu_{N+1}(\{(t,x) & \in K_-:\; w(t,x)>c(u)+\lambda\})\nonumber\\
\le &\;\mu_{N+1}(\{(t,x)\in K_-:
w(t,x)>c(u)+\lambda\;\,\mbox{and}\,\;W(t)\le c(u)+\lambda/2 \})\nonumber\\
&\; +\mu_{N+1}(\{(t,x)\in K_-:\,W(t)> c(u)+\lambda/2
\})=:I_1+I_2,\quad \lambda>0,\label{mainleft}\\
\mu_{N+1}(\{(t,x) & \in K_+:\; w(t,x)<c(u)-\lambda\})\nonumber\\
\le &\;\mu_{N+1}(\{(t,x)\in K_+:
w(t,x)<c(u)-\lambda\;\,\mbox{and}\,\;W(t)\ge c(u)-\lambda/2 \})\nonumber\\
&\; +\mu_{N+1}(\{(t,x)\in K_+:\,W(t)< c(u)-\lambda/2
\})=:I_3+I_4,\quad \lambda>0,\label{mainright}
\end{align}
and to estimate each of the four terms $I_j$ separately.

We begin with the estimates for $W$. To estimate $I_2$ and $I_4$ we
adopt some of the ideas developed in \cite{Za}. We set $J_-:=(0,\eta
\tau r^{2/\alpha})$, $J_+:=(\eta\tau r^{2/\alpha},\tau
r^{2/\alpha})$, and introduce for $\lambda>0$ the sets
$J_-(\lambda):=\{t\in J_-:\,W(t)> c(u)+\lambda \}$ and
$J_+(\lambda):=\{t\in J_+:\,W(t)< c(u)-\lambda \}$.

Interestingly, positivity and integrability of the function $e^W$
are sufficient to derive the desired estimate for $I_2$, cf. also
\cite[Theorem 2.3]{Za}. In fact, with $\rho=\tau r^{2/\alpha}$ we
have
\begin{align*}
e^\lambda \mu_1\big(J_-(\lambda)\big) & =
e^\lambda\mu_1\big(\{t\in J_-:\,
e^{W(t)}>e^{c(u)}e^{\lambda}\}\big)=\int_{J_-(\lambda)}e^\lambda \,dt\\
& \le \int_{J_-(\lambda)}e^{W(t)-c(u)} \,dt\le
\int_{J_-}e^{W(t)-c(u)} \,dt\\
& = \,\frac{g_{2-\alpha}(\eta\rho)}{(g_{1-\alpha}\ast
e^W)(\eta\rho)}\,\int_0^{\eta\rho}e^{W(t)}\,dt\\
& \le \,\frac{g_{2-\alpha}(\eta\rho)}{(g_{1-\alpha}\ast
e^W)(\eta\rho)}\,\cdot\,\frac{1}{g_{1-\alpha}(\eta\rho)}\,
\int_0^{\eta\rho}g_{1-\alpha}(\eta\rho-t)e^{W(t)}\,dt\\
& =
\,\frac{\Gamma(1-\alpha)}{\Gamma(2-\alpha)}\,\eta\rho=\,\frac{\eta\tau
r^{2/\alpha}}{1-\alpha},
\end{align*}
and therefore
\begin{equation} \label{I2est}
I_2=\mu_1\big(J_-(\lambda/2)\big)\mu_N(\delta B)\le\,\frac{2\eta\tau
\delta^N}{(1-\alpha)\lambda}\,r^{2/\alpha}\mu_N(B),\quad \lambda>0.
\end{equation}

We come now to $I_4$. For $m>0$ define the function $H_m$ on $\iR$
by $H_m(y)=y$, $y\le m$, and $H_m(y)=m+(y-m)/(y-m+1)$, $y\ge m$.
Then $H_m$ is increasing, concave, and bounded above by $m+1$.
Further, we have $H_m\in C^1(\iR)$, and so by concavity
\begin{equation} \label{Hmprop}
0\le yH_m'(y)\le H_m(y)\le m+1,\quad y\ge 0.
\end{equation}
Multiplying (\ref{log7}) by $e^W H_m'\big(e^W\big)$ and employing
(\ref{Hmprop}) as well as the fundamental identity
(\ref{fundidentity}), we infer that
\begin{equation} \label{West1}
\partial_t\Big(g_{1-\alpha,n}\ast
H_m\big(e^W\big)\Big)+\,\frac{C_1}{r^2}\,H_m\big(e^W\big)\ge - S_n
e^W H_m'\big(e^W\big),\quad \mbox{a.a.}\;t\in J.
\end{equation}
For $t\in J_+$ we shift the time by setting $s=t-\eta\tau
r^{2/\alpha}=t-\eta \rho$ and put $\tilde{f}(s)=f(s+\eta \rho)$,
$s\in (0,(1-\eta)\rho)$, for functions $f$ defined on $J_+$. By
the time-shifting identity (\ref{shiftprop}), (\ref{West1})
implies that for a.a. $s\in (0,(1-\eta)\rho)$
\begin{equation} \label{West2}
\partial_s\Big(g_{1-\alpha,n}\ast
H_m\big(e^{\tilde{W}}\big)\Big)+\,\frac{C_1}{r^2}\,H_m\big(e^{\tilde{W}}\big)\ge
\Upsilon_{n,m}(s)- \tilde{S}_n e^{\tilde{W}}
H_m'\big(e^{\tilde{W}}\big),
\end{equation}
with the history term
\[
\Upsilon_{n,m}(s)=\int_0^{\eta\rho}\big[-\dot{g}_{1-\alpha,n}(s+\eta\rho-\sigma)\big]
H_m\big(e^{W(\sigma)}\big)\,d\sigma.
\]

For $\theta\ge 0$ define the kernel $r_{\alpha,\theta}\in
L_{1,loc}(\iR_+)$ by means of
\[ r_{\alpha,\,\theta}(t)+\theta (r_{\alpha,\,\theta}\ast g_\alpha)(t)=
g_\alpha(t),\quad t>0.\] Observe that $r_{\alpha,\,0}=g_\alpha$.
Since $g_\alpha$ is completely monotone, $r_{\alpha,\,\theta}$
enjoys the same property (cf. \cite[Chap. 5]{GLS}), in particular
$r_{\alpha,\,\theta}(s)>0$ for all $s>0$. Moreover, we have (see
e.g. \cite{Za})
\[
r_{\alpha,\,\theta}(s) =\Gamma(\alpha)
g_\alpha(s)\,E_{\alpha,\alpha}(-\theta s^\alpha),\quad s>0,
\]
where $E_{\alpha,\beta}$ denotes the generalized
Mittag-Leffler-function defined by
\[ E_{\alpha,\beta}(z)=\sum_{n=0}^\infty \;\frac{z^n}
{\Gamma(n\alpha+\beta)}\;,\quad z\in \iC.
\]

We put $\theta=C_1/r^2$ and convolve (\ref{West2}) with
$r_{\alpha,\,\theta}$. We have a.e. in $(0,(1-\eta)\rho)$
\begin{align*}
r_{\alpha,\,\theta}\ast \partial_s\Big(g_{1-\alpha,n} & \ast
H_m\big(e^{\tilde{W}}\big)\Big)\,=\partial_s\Big(r_{\alpha,\,\theta}\ast
g_{1-\alpha,n}\ast H_m\big(e^{\tilde{W}}\big)\Big)\\
&\,=\partial_s\Big([g_\alpha-\theta (r_{\alpha,\,\theta}\ast
g_\alpha)]\ast g_{1-\alpha,n}\ast H_m\big(e^{\tilde{W}}\big)\Big)\\
&\,=h_n\ast H_m\big(e^{\tilde{W}}\big) -\theta
r_{\alpha,\,\theta}\ast h_n\ast H_m\big(e^{\tilde{W}}\big),
\end{align*}
and so we obtain a.e. in $(0,(1-\eta)\rho)$
\begin{align} \label{West3}
h_n\ast H_m\big(e^{\tilde{W}}\big)\ge &
\,\,r_{\alpha,\,\theta}\ast \Upsilon_{n,m}-
r_{\alpha,\,\theta}\ast\big[\tilde{S}_n
e^{\tilde{W}} H_m'\big(e^{\tilde{W}}\big)\big]\nonumber\\
& \,\,+\theta h_n\ast r_{\alpha,\,\theta}\ast
H_m\big(e^{\tilde{W}}\big) -\theta r_{\alpha,\,\theta}\ast
H_m\big(e^{\tilde{W}}\big).
\end{align}
Sending $n\to \infty$ and selecting an appropriate subsequence, if
necessary, it follows that
\begin{equation} \label{West4}
H_m\big(e^{\tilde{W}}\big)\ge r_{\alpha,\,\theta}\ast
\Upsilon_{m},\quad\quad \mbox{a.a.}\;s\in(0,(1-\eta)\rho),
\end{equation}
where
\[
\Upsilon_{m}(s)=\int_0^{\eta\rho}\big[-\dot{g}_{1-\alpha}(s+\eta\rho-\sigma)\big]
H_m\big(e^{W(\sigma)}\big)\,d\sigma.
\]

Observe that for $s\in (0,(1-\eta)\rho)$ we have
\[
0\le \theta s^\alpha\le \,\frac{C_1}{r^2}\,(1-\eta)^\alpha \big(\tau
r^{2/\alpha}\big)^\alpha=C_1(1-\eta)^\alpha\tau^\alpha=:\omega,
\]
and thus by continuity and strict positivity of $E_{\alpha,\alpha}$
in $(-\infty,0]$,
\[
r_{\alpha,\,\theta}(s)\ge  \Gamma(\alpha) g_\alpha(s)\min_{y\in
[0,\omega]} E_{\alpha,\alpha}(-y)=:C_2(\alpha,\omega)\Gamma(\alpha)
g_\alpha(s),\quad s\in (0,(1-\eta)\rho).
\]
We may then argue as in \cite[Section 2.1]{Za} to obtain
\[
H_m\big(e^{\tilde{W}(s)}\big)\ge C_2(\alpha,\omega)\,\frac{\alpha
(s/[\eta\rho])^\alpha}{1+(s/[\eta\rho])}\,(\eta\rho)^{\alpha-1}\big(g_{1-\alpha}\ast
H_m\big(e^{W}\big)\big) (\eta\rho),\quad \mbox{a.a.}\;s\in
(0,(1-\eta)\rho).
\]
Evidently, $H_m(y)\nearrow y$ as $m\to \infty$ for all $y\in \iR$.
Thus by sending $m\to \infty$ and applying Fatou's lemma we conclude
that
\begin{equation} \label{West5}
e^{\tilde{W}(s)}\ge C_2(\alpha,\omega)\,\frac{\alpha
(s/[\eta\rho])^\alpha}{1+(s/[\eta\rho])}\,(\eta\rho)^{\alpha-1}\big(g_{1-\alpha}\ast
e^{W}\big) (\eta\rho),\quad \mbox{a.a.}\;s\in (0,(1-\eta)\rho).
\end{equation}

We then employ (\ref{West5}) to estimate as follows.
\begin{align*}
e^\lambda \mu_1\big(J_+(\lambda)\big) & = e^\lambda\mu_1\big(\{t\in
J_+:\,
e^{W(t)}<e^{c(u)}e^{-\lambda}\}\big)=\int_{J_+(\lambda)}e^\lambda \,dt\\
& \le \int_{J_+(\lambda)}e^{c(u)-W(t)} \,dt\le
\int_{J_+}e^{c(u)-W(t)} \,dt\\
& = \,\frac{(g_{1-\alpha}\ast
e^W)(\eta\rho)}{g_{2-\alpha}(\eta\rho)}\,\int_0^{(1-\eta)\rho}
e^{-\tilde{W}(s)}\,ds\\
& \le \,\frac{C_2(\alpha,\omega)^{-1}(\eta\rho)^{1-\alpha}}{\alpha
g_{2-\alpha}(\eta\rho)}\,\int_0^{(1-\eta)\rho}(1+s/\eta\rho)(s/\eta\rho)^{-\alpha}\,ds\\
& = \,\frac{\Gamma(2-\alpha)\eta\rho}{\alpha
C_2(\alpha,\omega)}\,\int_0^{\frac{1-\eta}{\eta}}
\sigma^{-\alpha}(1+\sigma)\,d\sigma=C_3(\alpha,\eta,\omega)\rho.
\end{align*}
Hence
\begin{equation} \label{I4est}
I_4=\mu_1\big(J_+(\lambda/2)\big)\mu_N(\delta
B)\le\,\frac{2C_3(\alpha,\eta,\omega)
\delta^N}{\lambda}\,r^{2/\alpha}\mu_N(B),\quad \lambda>0.
\end{equation}

We come now to $I_1$. Set $J_1(\lambda)=\{t\in
J_-:\,c-W(t)+\lambda/2\ge 0\}$ and $\Omega^-_t(\lambda)=\{x\in
\delta B:\,w(t,x)>c+\lambda\},\,t\in J_1(\lambda)$, where $c=c(u)$
is given by (\ref{cwahl}). For $t\in J_1(\lambda)$, we have
\[ w(t,x)-W(t)>c-W(t)+\lambda\ge \lambda/2,\quad x\in
\Omega^-_t(\lambda),\] and thus we deduce from (\ref{log7}) that
a.e. in $J_1(\lambda)$
\begin{equation} \label{log8}
\frac{\nu}{16r^2
\mu_N(B)}\,\,\mu_N\big(\Omega^-_t(\lambda)\big)\le
\frac{1}{(c-W+\lambda)^2}\,\Big(e^{-W}
\partial_t(g_{1-\alpha,n}\ast
e^W)+\,\frac{C_1}{r^2}\,+S_n\Big).
\end{equation}
Set $\chi(t,\lambda)=\mu_N\big(\Omega^-_t(\lambda)\big)$, if $t\in
J_1(\lambda)$, and $\chi(t,\lambda)=0$ in case $t\in J_-\setminus
J_1(\lambda)$. Let further $H(y)=(c-\log y+\lambda)^{-1},\,0<y\le
y_*:=e^{c+\lambda/2}$. Clearly, $H'(y)= (c-\log
y+\lambda)^{-2}y^{-1}$ as well as
\[ H''(y)=\,\frac{1}{(c-\log y+\lambda)^2 y^2}\,\Big(\frac{2}{c-\log
y+\lambda}-1\Big),\quad 0<y\le y_*,\] which shows that $H$ is
concave in $(0,y_*]$ whenever $\lambda\ge 4$. We will assume this in
what follows.

We next choose a $C^1$ extension $\bar{H}$ of $H$ on $(0,\infty)$
such that $\bar{H}$ is concave, $0\le\bar{H}'(y)\le
\bar{H}'(y_*),\,y_*\le y \le 2 y_*$, and $\bar{H}'(y)=0,\,y\ge 2
y_*$. Then
\begin{equation} \label{log8a}
 0\le y\bar{H}'(y)\le
\,\frac{2}{\lambda},\quad y>0.
\end{equation}
In fact, for $y\in(0,y_*]$ we have
\begin{equation} \label{log8aa}
y\bar{H}'(y)=\,\frac{1}{(c-\log y+\lambda)^2}\,\le
\,\frac{1}{(c-\log y_*+\lambda)^2}\,\le \,\frac{4}{\lambda^2}\,\le
\,\frac{1}{\lambda},
\end{equation}
while in case $y\in[y_*,2y_*]$ we may simply estimate
\[
y\bar{H}'(y)\le 2y_*\bar{H}'(y_*)\le \,\,\frac{2}{\lambda}.
\]

It is clear that $\bar{H}$ is bounded above. There holds
\begin{equation} \label{log8b}
\bar{H}(y)\le \,\frac{3}{\lambda},\quad y>0.
\end{equation}
To see this, note that since $\bar{H}$ is nondecreasing with
$\bar{H}'(y)=0$ for all $y\ge 2 y_*$, the claim follows if the
inequality is valid for all $y\in [y_*,2y_*]$. For such $y$ we
have by (\ref{log8aa}) and by concavity of $\bar{H}$
\[
\bar{H}(y)\le \bar{H}(y_*)+\bar{H}'(y_*)(y-y_*)\le
\bar{H}(y_*)+y_*\bar{H}'(y_*)\le \,\frac{3}{\lambda}.
\]

Observe also that
\[ e^{W(t)} H'(e^{W(t)})=\,\frac{1}{(c-W(t)+\lambda)^2
},\quad \mbox{a.a.}\;t\in J_1(\lambda).\] Since $\bar{H}'\ge 0$,
and $e^{-W}
\partial_t(g_{1-\alpha,n}\ast
e^W)+C_1 r^{-2}+S_n\ge 0$ on $J_-$ by virtue of (\ref{log7}), we
infer from (\ref{log8}) and (\ref{log8a}) that
\begin{align} \label{log9}
\frac{\nu}{16r^2 \mu_N(B)}\,\,\chi(t,\lambda) & \le
e^W\bar{H}'(e^{W})\Big(e^{-W}\partial_t(g_{1-\alpha,n}\ast
e^W)+\,\frac{C_1}{r^2}\,+S_n\Big)\nonumber\\
& \le \bar{H}'(e^{W})\partial_t(g_{1-\alpha,n}\ast
e^W)+\,\frac{2C_1}{\lambda r^2}\,+\,\frac{2|S_n(t)|}{\lambda},
\quad\mbox{a.a.}\; t\in J_-.
\end{align}
Since $\bar{H}$ is concave, the fundamental identity
(\ref{fundidentity}) yields
\begin{align*}
\bar{H}'(e^{W})\partial_t(g_{1-\alpha,n}\ast e^W) & \le
\partial_t\Big(g_{1-\alpha,n}\ast \bar{H}\big(e^W\big)\Big)
+ \Big(-\bar{H}(e^{W})+\bar{H}'(e^{W})
e^W\Big)g_{1-\alpha,n}\\
& \le \partial_t\Big(g_{1-\alpha,n}\ast
\bar{H}\big(e^W\big)\Big)+\,\frac{2}{\lambda}\,g_{1-\alpha,n},\quad
\mbox{a.a.}\; t\in J_-,
\end{align*}
which, together with (\ref{log9}), gives a.e. in $J_-$
\begin{equation} \label{log10}
\frac{\nu}{16r^2 \mu_N(B)}\,\,\chi(t,\lambda) \le
\partial_t\Big((g_{1-\alpha,n}\ast \bar{H}\big(e^W\big)\Big)
+\,\frac{2}{\lambda}\,g_{1-\alpha,n}+ \,\frac{2C_1}{\lambda
r^2}\,+\,\frac{2|S_n(t)|}{\lambda}\,.
\end{equation}
We then integrate (\ref{log10}) over $J_-=(0,\eta \rho)$ and
employ (\ref{log8b}) for the estimate
\[ \Big(g_{1-\alpha,n}\ast
\bar{H}\big(e^W\big)\Big)(\eta\rho)\le
\,\frac{3}{\lambda}\,\int_0^{\eta\rho}g_{1-\alpha,n}(t)\,dt.
\]
By sending $n\to \infty$, this leads to
\begin{align*}
\int_{J_1(\lambda)}\mu_N  & \big(\Omega^-_t(\lambda)\big) \,dt
  = \int_0^{\eta\rho} \chi(t,\lambda)\,dt
 \le\,\frac{16r^2\mu_N(B)}{\nu}\,\Big(\frac{5}{\lambda}\,g_{2-\alpha}(\eta\rho)+
 \,\frac{2C_1\eta\rho}{\lambda
r^2}\Big)\\
& =\,\frac{16r^{2/\alpha}\mu_N(B)}{\nu\lambda}\,
\big(5g_{2-\alpha}(\eta\tau)+
2C_1\eta\tau\big)=:C_4\,\frac{r^{2/\alpha}\mu_N(B)}{\lambda},\quad\lambda\ge
4.
\end{align*}
Hence with $C_5=\max\{4\tau,C_4\}$ we find that
\begin{equation} \label{I1est}
I_1\le \,\frac{C_5 r^{2/\alpha}\mu_N(B)}{\lambda},\quad\lambda>0.
\end{equation}

It remains to derive the desired estimate for $I_3$. To this
purpose we shift again the time by putting $s=t-\eta\rho$, and
denote the corresponding transformed functions as above by
$\tilde{W}$, $\tilde{w}$, ... and so forth. Set further
$\tilde{J}_+:=(0,(1-\eta)\rho)$. By the time-shifting property
(\ref{shiftprop}) and by positivity of $e^W$, relation
(\ref{log7}) then implies
\begin{equation} \label{log10a}
\frac{\nu}{16r^2 \mu_N(B)}\,\int_{\delta B}
(\tilde{w}-\tilde{W})^2 dx \le e^{-\tilde{W}}
\partial_s(g_{1-\alpha,n}\ast
e^{\tilde{W}})+\,\frac{C_1}{r^2}\,+\tilde{S}_n(s),\quad
\mbox{a.a.}\;s\in \tilde{J}_+.
\end{equation}
Next, set $J_2(\lambda)=\{s\in
\tilde{J}_+:\tilde{W}(s)-c+\lambda/2\ge 0\}$ and
$\Omega_{s}^+(\lambda)=\{x\in \delta B:
\tilde{w}(s,x)<c-\lambda\},\,s\in J_2(\lambda)$. For $s\in
J_2(\lambda)$, we have
\[ \tilde{W}(s)-\tilde{w}(s,x)\ge
\tilde{W}(s)-c+\lambda\ge \lambda/2,\quad
x\in\Omega_{s}^+(\lambda),
\]
and thus (\ref{log10a}) yields that a.e. in $J_2(\lambda)$
\begin{equation} \label{log12}
\frac{\nu}{16r^2
\mu_N(B)}\,\,\mu_N\big(\Omega^+_s(\lambda)\big)\le
\frac{1}{(\tilde{W}-c+\lambda)^2}\,\Big(e^{-\tilde{W}}
\partial_s(g_{1-\alpha,n}\ast
e^{\tilde{W}})+\,\frac{C_1}{r^2}\,+\tilde{S_n}\Big).
\end{equation}

We proceed now similarly as above for the term $I_1$. Set
$\chi(s,\lambda)=\mu_N\big(\Omega^+_{s}(\lambda)\big)$, if $s\in
J_2(\lambda)$, and $\chi(s,\lambda)=0$ in case $s\in
\tilde{J}_+\setminus J_1(\lambda)$. We consider this time the
convex function $H(y)=(\log y-c+\lambda)^{-1}$ for $y\ge
y_*:=e^{c-\lambda/2}$ with derivative $H'(y)=-(\log
y-c+\lambda)^{-2} y^{-1}<0$. We define a $C^1$ extension $\bar{H}$
of $H$ on $[0,\infty)$ by means of
\[
\bar{H}(y)=\left\{ \begin{array}{l@{\;:\;}l}
H'(y_*)(y-y_*)+H(y_*) & 0\le y< y_* \\
H(y) & y\ge y_*.
\end{array} \right.
\]
Evidently, $-\bar{H}$ is concave in $[0,\infty)$ and
\begin{equation} \label{log13}
0\le-\bar{H}'(y)y \le \,\frac{1}{(\log
y_*-c+\lambda)^2}\,\le\,\frac{1}{(\lambda/2)^2}\,\le
\,\frac{4}{\lambda},\quad y \ge 0,\;\lambda\ge 1.
\end{equation}
We will assume $\lambda\ge 1$ in the subsequent lines.

Observe that
\[ -e^{\tilde{W}(s)} H'(e^{\tilde{W}(s)})=\,\frac{1}{(\tilde{W}(s)-c+\lambda)^2
},\quad \mbox{a.a.}\;s\in J_2(\lambda).\]
Since $-\bar{H}'\ge 0$,
and $e^{-\tilde{W}}
\partial_s(g_{1-\alpha,n}\ast
e^{\tilde{W}})+C_1 r^{-2}+\tilde{S_n}\ge 0$ on $\tilde{J}_+$ due
to (\ref{log10a}), it thus follows from (\ref{log12}) and
(\ref{log13}) that
\begin{align} \label{log14}
\frac{\nu}{16r^2 \mu_N(B)}\,\,\chi(s,\lambda) & \le
-e^{\tilde{W}}\bar{H}'(e^{\tilde{W}})\Big(e^{-\tilde{W}}\partial_s(g_{1-\alpha,n}\ast
e^{\tilde{W}})+\,\frac{C_1}{r^2}\,+\tilde{S}_n\Big)\nonumber\\
& \le -\bar{H}'(e^{\tilde{W}})\partial_s(g_{1-\alpha,n}\ast
e^{\tilde{W}})+\,\frac{4C_1}{\lambda
r^2}\,+\,\frac{4|\tilde{S}_n(s)|}{\lambda}, \quad\mbox{a.a.}\;
s\in \tilde{J}_+.
\end{align}
By concavity of $-\bar{H}$, the fundamental identity
(\ref{fundidentity}) provides the estimate
\begin{align*}
-\bar{H}'(e^{\tilde{W}}) & \partial_s(g_{1-\alpha,n}\ast
e^{\tilde{W}})  \le -\partial_s\Big(g_{1-\alpha,n}\ast
\bar{H}\big(e^{\tilde{W}}\big)\Big)+\Big(\bar{H}(e^{\tilde{W}})-
\bar{H}'(e^{\tilde{W}})e^{\tilde{W}}\Big)g_{1-\alpha,n} \\
& \le -\partial_s\Big(g_{1-\alpha,n}\ast
\bar{H}\big(e^{\tilde{W}}\big)\Big)+\bar{H}(0)g_{1-\alpha,n} \le
-\partial_s\Big(g_{1-\alpha,n}\ast
\bar{H}\big(e^{\tilde{W}}\big)\Big)
+\,\frac{6}{\lambda}\,g_{1-\alpha,n},
\end{align*}
a.e. in $\tilde{J}_+$, which when combined with (\ref{log14})
leads to
\[
\frac{\nu}{16r^2 \mu_N(B)}\,\,\chi(s,\lambda) \le
-\partial_s\Big(g_{1-\alpha,n}\ast
\bar{H}\big(e^{\tilde{W}}\big)\Big)
+\,\frac{6}{\lambda}\,g_{1-\alpha,n}+\,\frac{4C_1}{\lambda
r^2}\,+\,\frac{4|\tilde{S}_n(s)|}{\lambda},
\]
for a.a. $s\in \tilde{J}_+$. We integrate this estimate over
$\tilde{J}_+$ and send $n\to \infty$ to the result
\begin{align*}
\int_{J_2(\lambda)}\mu_N\big(  & \Omega^+_s(\lambda)\big) \,ds
  = \int_0^{(1-\eta)\rho} \!\!\!\!\chi(s,\lambda)\,ds\le
  \,\frac{16r^2\mu_N(B)}{\nu}\,\Big(\,\frac{6}{\lambda}\,g_{2-\alpha}\big((1-\eta)\rho\big)+\,\frac{4C_1(1-\eta)\rho}{\lambda
  r^2}\,\Big)\\
  &
  =\,\frac{16r^{2/\alpha}\mu_N(B)}{\nu\lambda}\big(6g_{2-\alpha}\big((1-\eta)\tau\big)+4C_1(1-\eta)\tau\big)=:C_6\,
  \frac{r^{2/\alpha}\mu_N(B)}{\lambda},\quad \lambda\ge 1.
\end{align*}
Hence with $C_7=\max\{\tau,C_6\}$ we obtain that
\begin{equation} \label{I3est}
I_3\le \,\frac{C_7 r^{2/\alpha}\mu_N(B)}{\lambda},\quad\lambda>0.
\end{equation}

Finally, combining (\ref{mainleft}), (\ref{mainright}), and
(\ref{I2est}), (\ref{I4est}), (\ref{I1est}), (\ref{I3est})
establishes the theorem. $\square$
\subsection{The final step}
We are now in position to prove Theorem \ref{localweakHarnack}.
Without loss of generality we may assume that $u\ge \varepsilon$ for
some $\varepsilon>0$; otherwise replace $u$ by $u+\varepsilon$,
which is a supersolution of (\ref{MProb}) with $u_0+\varepsilon$
instead of $u_0$, and eventually let $\varepsilon\to 0+$.

For $0<\sigma\le 1$, we set $U_\sigma=(t_0+(2-\sigma)\tau
r^{2/\alpha},t_0+2\tau r^{2/\alpha})\times \sigma B$ and
$U'_\sigma=(t_0,t_0+\sigma\tau r^{2/\alpha})\times \sigma B$.
Clearly, $Q_-(t_0,x_0,r)=U'_\delta$ and $Q_+(t_0,x_0,r)=U_\delta$.

By Theorem \ref{superest1},
\[
\esup_{U_{\sigma'}}{u^{-1}} \le \Big(\frac{C \mu_{N+1}(U_1)^{-1}
}{(\sigma-\sigma')^{\tau_0}}\Big)^{1/\gamma}
|u^{-1}|_{L_{\gamma}(U_\sigma)},\quad \delta\le \sigma'<\sigma\le
1,\; \gamma\in (0,1].
\]
Here $C=C(\nu,\Lambda,\delta,\tau,\alpha,N)$ and
$\tau_0=\tau_0(\alpha,N)$. This shows that the first hypothesis of
Lemma \ref{abslemma} is satisfied by any positive constant multiple
of $u^{-1}$ with $\beta_0=\infty$.

Consider now $f_1=u^{-1}e^{c(u)}$ where $c(u)$ is the constant from
Theorem \ref{logest} with $K_-=U'_1$ and $K_+=U_1$. Since $\log
f_1=c(u)-\log u$, we see from Theorem \ref{logest}, estimate
(\ref{logestright}), that
\[ \mu_{N+1}(\{(t,x)\in U_1:\;\log f_1(t,x)>\lambda\})\le
M\mu_{N+1}(U_1)\lambda^{-1},\quad \lambda>0,\] where
$M=M(\nu,\Lambda,\delta,\tau,\eta,\alpha,N)$. Hence we may apply
Lemma \ref{abslemma} with $\beta_0=\infty$ to $f_1$ and the family
$U_\sigma$; thereby we obtain
\[
\esup_{U_\delta} f_1\le M_1\] with
$M_1=M_1(\nu,\Lambda,\delta,\tau,\eta,\alpha,N)$. In terms of $u$
this means that
\begin{equation} \label{HH1}
e^{c(u)}\le M_1\, \einf_{U_\delta} u.
\end{equation}

On the other hand, Theorem \ref{superest2} yields
\[
|u|_{L_{p}(U_{\sigma'}')}\le \Big(\frac{C\mu_{N+1}(U'_1)^{-1}
}{(\sigma-\sigma')^{\tau_1}}\Big)^{1/\gamma-1/p}
|u|_{L_{\gamma}(U'_\sigma)},\quad \delta\le \sigma'<\sigma\le 1,\;
0<\gamma\le p/\tilde{\kappa}.
\]
Here $C=C(\nu,\Lambda,\delta,\tau,\alpha,N,p)$ and
$\tau_1=\tau_1(\alpha,N)$. Thus the first hypothesis of Lemma
\ref{abslemma} is satisfied by any positive constant multiple of $u$
with $\beta_0=p$ and $\eta=1/\tilde{\kappa}$. Taking $f_2=u
e^{-c(u)}$ with $c(u)$ from above, we have $\log f_2=\log u-c(u)$
and so Theorem \ref{logest}, estimate (\ref{logestleft}), gives
\[ \mu_{N+1}(\{(t,x)\in U'_1:\;\log f_2(t,x)>\lambda\})\le
M\mu_{N+1}(U'_1)\lambda^{-1},\quad \lambda>0,\] where $M$ is as
above. Therefore we may again apply Lemma \ref{abslemma}, this time
to the function $f_2$ and the sets $U'_\sigma$, and with $\beta_0=p$
and $\eta=1/\tilde{\kappa}$; we get
\[
|f_2|_{L_p(U'_\delta)}\le M_2 \mu_{N+1}(U'_1)^{1/p},
\]
where $M_2=M_2(\nu,\Lambda,\delta,\tau,\eta,\alpha,N,p)$. Rephrasing
then yields
\begin{equation} \label{HH2}
\mu_{N+1}(U'_1)^{-1/p}|u|_{L_p(U'_\delta)}\le M_2 e^{c(u)}.
\end{equation}

Finally, we combine (\ref{HH1}) and (\ref{HH2}) to the result
\[ \mu_{N+1}(U'_1)^{-1/p}|u|_{L_p(U'_\delta)}\le M_1 M_2\, \einf_{U_\delta}
u,
\]
which proves the assertion. $\square$
\section{Optimality of the exponent
$\frac{2+N\alpha}{2+N\alpha-2\alpha}$ in the weak Harnack
inequality}
In this section we will show that the exponent
$\frac{2+N\alpha}{2+N\alpha-2\alpha}$ in Theorem
\ref{localweakHarnack} is optimal.

To this purpose consider the nonhomogeneous fractional diffusion
equation on $\iR^N$
\begin{equation} \label{opt1}
\partial_t^\alpha u-\Delta u=f,\quad t\in (0,T],\,x\in \iR^N,
\end{equation}
with initial condition
\begin{equation} \label{opt2}
u(0,x)=0,\quad x\in \iR^n.
\end{equation}
Following \cite{Koch}, we say that a function $u\in C([0,T]\times
\iR^N)\cap C((0,T];C^2(\iR^N))$ with $g_{1-\alpha}\ast u\in
C^1((0,T];C(\iR^N))$ is a classical solution of the problem
(\ref{opt1}), (\ref{opt2}) if $u$ satisfies (\ref{opt1}) and
(\ref{opt2}). For any bounded continuous function $f$ that is
locally H\"older continuous in $x$, there exists a unique classical
solution $u$ of the problem (\ref{opt1}), (\ref{opt2}), and it is of
the form
\begin{equation} \label{opt3}
u(t,x)=\int_0^t \int_{\iR^N} Y(t-\tau,x-y)f(\tau,y)\,dy\,d\tau,
\end{equation}
where
\[
Y(t,x)=c(N)|x|^{-N}t^{\alpha-1}
H^{20}_{12}\Big(\frac{1}{4}\,t^{-\alpha}|x|^2\Big|{}^{(\alpha,\alpha)}_{(N/2,1),\,(1,1)}\Big),
\]
cf. \cite{Koch}. Here
$H^{20}_{12}(z|{}^{(\alpha,\alpha)}_{(N/2,1),\,(1,1)})$ denotes a
special $H$ function (also termed Fox's $H$ function), see
\cite[Section 1.12]{KST} and \cite{Koch} for its definition. It is
differentiable for $z>0$, the asymptotic behaviour for $z\to \infty$
and $z\to +0$, respectively, is described in \cite[formulae (3.9)
and (3.14)]{Koch}. It has been also proved in \cite{Koch} that $Y$
is nonnegative.

We choose a smooth and nonnegative approximation of unity
$\{\phi_n(t,x)\}_{n\in \iN}$ in $\iR_+\times \iR^N$ such that each
$\phi_n$ is bounded. Put $f=\phi_n$ in (\ref{opt1}) and denote the
corresponding classical solution of (\ref{opt1}), (\ref{opt2}) by
$u_n$. Evidently, $u_n$ is nonnegative and satisfies
\[
\partial_t^\alpha u_n-\Delta u_n=\phi_n\ge  0,\quad t\in (0,T],\,x\in
\iR^N.
\]
Hence $u_n$ is a nonnegative supersolution of (\ref{opt1}) with
$f=0$ for all $n\in \iN$.

Suppose the weak Harnack inequality (\ref{localwHarnackF}) holds for
some $p\ge \frac{2+N\alpha}{2+N\alpha-2\alpha}$. Then, by taking
$Q_-=(0,1)\times B(0,1)$ and $Q_+=(2,3)\times B(0,1)$ it follows
that
\begin{equation} \label{opt4}
\big(\int_{Q_-}u_n^p\,d\mu_{N+1}\big)^{1/p}\le C \inf_{Q_+}
u_n,\quad n\in \iN,
\end{equation}
where the constant $C$ is independent of $n$. Since $u_n\to Y$ in
the distributional sense as $n\to \infty$, we have
\[
\inf_{Q_+} u_n\le \,\frac{1}{\mu_{N+1}(Q_+)}\,\int_{Q_+}u_n
\,d\mu_{N+1}\le 1+ \,\frac{1}{\mu_{N+1}(Q_+)}\,\int_{Q_+}Y
\,d\mu_{N+1}<\infty,\quad n\ge n_0,
\]
for a sufficiently large $n_0$. On the other hand, the left-hand
side of (\ref{opt4}) cannot stay bounded, since $Y\notin L_p(Q_-)$
for $p\ge \frac{2+N\alpha}{2+N\alpha-2\alpha}$. In fact, writing
$H^{20}_{12}(z)=H^{20}_{12}(z|{}^{(\alpha,\alpha)}_{(N/2,1),\,(1,1)})$
for short, we have
\begin{align*}
|Y|_{L_p(Q_-)}^p & = \int_0^1 \int_{B(0,1)} c(N)^p
|x|^{-Np}t^{(\alpha-1)p}
H^{20}_{12}\big(t^{-\alpha}|x|^2/4\big)^p\,dx\,dt\\
& = c_1 \int_0^1 \int_0^1
r^{N-1-Np}t^{(\alpha-1)p}H^{20}_{12}\big(t^{-\alpha}r^2/4\big)^p\,dr\,dt\\
& = c_1 \int_0^1 \int_0^{t^{-\alpha/2}}\big(\rho
t^{\alpha/2})^{N-1-Np}t^{(\alpha-1)p+\alpha/2}H^{20}_{12}\big(\rho^2/4\big)^p\,d\rho\,dt\\
& \ge c_1 \int_0^1 t^{\alpha(N-Np)/2+(\alpha-1)p}\,dt \,\int_0^1
\rho^{N-1-Np}H^{20}_{12}\big(\rho^2/4\big)^p\,d\rho\\
& \ge c_2 \int_0^1 t^{\alpha(N-Np)/2+(\alpha-1)p}\,dt,
\end{align*}
with some positive constant $c_2$. The last integral diverges for
all $p\ge \frac{2+N\alpha}{2+N\alpha-2\alpha}$. Hence (\ref{opt4})
yields a contradiction.
\section{Applications of the weak Harnack inequality}
The strong maximum principle for weak subsolutions of (\ref{MProb})
may be easily derived as a consequence of the weak Harnack
inequality.
\begin{satz} \label{strongmax}
Let $\alpha\in(0,1)$, $T>0$, and $\Omega\subset \iR^N$ be a bounded
domain. Suppose the assumptions (H1)--(H3) are satisfied. Let $u\in
Z_\alpha$ be a weak subsolution of (\ref{MProb}) in $\Omega_T$ and
assume that $0\le \esup_{\Omega_T}u<\infty$ and that $\esup_{\Omega}
u_0\le \esup_{\Omega_T}u$. Then, if for some cylinder
$Q=(t_0,t_0+\tau r^{2/\alpha})\times B(x_0,r)\subset \Omega_T$ with
$t_0,\tau,r>0$ and $\overline{B(x_0,r)}\subset \Omega$ we have
\begin{equation} \label{strrel}
\esup_{Q}u \,=\,\esup_{\Omega_T}u,
\end{equation}
the function $u$ is constant on $(0,t_0)\times \Omega$.
\end{satz}
{\em Proof:} Let $M=\esup_{\Omega_T}u$. Then $v:=M-u$ is a
nonnegative weak supersolution of (\ref{MProb}) with $u_0$ replaced
by $v_0:=M-u_0\ge 0$. For any $0\le t_1< t_1+\eta r^{2/\alpha}<t_0$
the weak Harnack inequality with $p=1$ applied to $v$ yields an
estimate of the form
\[
r^{-(N+2/\alpha)}\int_{t_1}^{t_1+\eta
r^{2/\alpha}}\int_{B(x_0,r)}(M-u)\,dx\,dt\le C\,\einf_Q (M-u)\,=\,0.
\]
This shows that $u=M$ a.e. in $(0,t_0)\times B(x_0,r)$. As in the
classical parabolic case (cf. \cite{Lm}) the assertion now follows
by a chaining argument. $\square$

$\mbox{}$

We next apply the weak Harnack inequality to establish continuity at
$t=0$ for weak solutions.
\begin{satz} \label{Hoeldert=0}
Let $\alpha\in(0,1)$, $T>0$, and $\Omega\subset \iR^N$ be a bounded
domain. Suppose the assumptions (H1) and (H2) are satisfied. Let
$u\in Z_\alpha$ be a bounded weak solution of (\ref{MProb}) in
$\Omega_T$ with $u_0=0$. Then $u$ is continuous at $(0,x_0)$ for all
$x_0\in \Omega$ and $\lim_{(t,x)\to (0,x_0)}u(t,x)=0$. Moreover,
letting $\eta>0$ we have for any cylinder $Q(x_0,r_0):=(0,\eta
r_0^{2/\alpha})\times B(x_0,r_0)\subset \Omega_T$ and $r\in (0,r_0]$
\begin{equation} \label{oscest}
\eosc_{Q(x_0,r)} u \le C\Big(\,\frac{r}{r_0}\,\Big)^\delta
|u|_{L_\infty(\Omega_T)},
\end{equation}
with $\eosc_{Q(x_0,r)}=\esup_{Q(x_0,r)}-\einf_{Q(x_0,r)}$ and
constants $C=C(\nu,\Lambda,\eta,\alpha,N)>0$ and
$\delta=\delta(\nu,\Lambda,\eta,\alpha,N)\in (0,1)$.
\end{satz}
{\em Proof:} Let $u\in Z_\alpha$ be a bounded weak solution of
(\ref{MProb}) in $\Omega_T$ with $u_0=0$. Set $u(t,x)=0$ and
$A(t,x)=Id$ for $t<0$ and $x\in \Omega$. For $T_0>0$ we shift the
time by setting $s=t+T_0$ and put $\tilde{f}(s)=f(s-T_0)$, $s\in
(0,T+T_0)$, for functions $f$ defined on $(-T_0,T)$. Since
$Du(t,\cdot)=0$ for $t<0$ and
\[
\partial_t(g_{1-\alpha,n}\ast u)(t,x)=\partial_t\int_{-T_0}^t
g_{1-\alpha,n}(t-\tau)u(\tau,x)\,d\tau=\partial_s(g_{1-\alpha,n}\ast
\tilde{u})(s,x),\] the function $\tilde{u}$ is a bounded weak
solution of
\[
\partial_s^\alpha \tilde{u}-\mbox{div}\,\big(\tilde{A}(s,x)D\tilde{u}\big)=0,\quad s\in (0,T+T_0),\,x\in
\Omega.
\]

Next, assuming $r\in (0,r_0/2]$ we introduce the cylinders
\begin{align*}
Q_*(x_0,r) & =\big(-\eta r^{2/\alpha},\eta r^{2/\alpha}\big)\times
B(x_0,r),\\
Q_-(x_0,r) & =\big(-\eta (2r)^{2/\alpha},-\eta
(3r/2)^{2/\alpha}\big)\times B(x_0,r),
\end{align*}
and denote by $\tilde{Q}_*(x_0,r)$ resp. $\tilde{Q}_-(x_0,r)$ the
corresponding cylinders in the $(s,x)$ coordinate system. Let us
write $M_i=\esup_{\tilde{Q}_*(x_0,ir)}\tilde{u}$ and
$m_i=\einf_{\tilde{Q}_*(x_0,ir)}\tilde{u}$ for $i=1,2$. Choosing
$T_0\ge\eta (2r)^{2/\alpha}$, we may apply Theorem
\ref{localweakHarnack} with $p=1$ to the functions $M_2-\tilde{u}$,
$\tilde{u}-m_2$, which are nonnegative in $(0,\eta
(2r)^{2/\alpha}+T_0)\times B(x_0,2r)$, thereby obtaining
\begin{align*}
r^{-N+2/\alpha}\int_{\tilde{Q}_-(x_0,r)}(M_2-\tilde{u})\,d\mu_{N+1}
& \le C(M_2-M_1),\\
r^{-N+2/\alpha}\int_{\tilde{Q}_-(x_0,r)}(\tilde{u}-m_2)\,d\mu_{N+1}
& \le C(m_1-m_2),
\end{align*}
where $C>1$ is a constant independent of $u$ and $r$. By addition,
it follows that
\[
M_2-m_2\le C(M_2-m_2+m_1-M_1).
\]
Writing
$\omega(x_0,r)=\esup_{\tilde{Q}_*(x_0,ir)}\tilde{u}-\einf_{\tilde{Q}_*(x_0,ir)}\tilde{u}$,
this yields
\begin{equation} \label{Hc1}
\omega(x_0,r)\le \theta \omega(x_0,2r),\quad r\le r_0/2,
\end{equation}
where $\theta=1-C^{-1}\in (0,1)$. Iterating (\ref{Hc1}) as in the
proof of \cite[Lemma 8.23]{GilTrud} we obtain
\[
\omega(x_0,r)\le \,\frac{1}{\theta}\,\Big(\frac{r}{r_0}\Big)^{\log
\theta/\log(1/2)}\omega(x_0,r_0),\quad r\le r_0.
\]
The estimate (\ref{oscest}) then follows by transforming back to the
function $u$ and using that $u=0$ for negative times. In particular,
we also see that $u$ is continuous at $(0,x_0)$ for all $x_0\in
\Omega$ and that $\lim_{(t,x)\to (0,x_0)}u(t,x)=0$. $\square$

$\mbox{}$

The last application is a theorem of Liouville type. We say that a
function $u$ on $\iR_+\times \iR^N$ is a {\em global weak solution}
of
\begin{equation} \label{GlE}
\partial_t^\alpha u-\mbox{div}\,\big(A(t,x)Du\big)=0,
\end{equation}
if it is a weak solution of (\ref{GlE}) in $(0,T)\times B(0,r)$ for
all $T>0$ and $r>0$.
\begin{korollar} \label{Liouville}
Let $\alpha\in(0,1)$. Assume that $A\in L_\infty(\iR_+\times
\iR^N;\iR^{N\times N})$ and that there exists $\nu>0$ such that
\[
\big(A(t,x)\xi|\xi\big)\ge \nu|\xi|^2,\quad\mbox{for a.a.}\;
(t,x)\in\iR_+\times \iR^N,\; \mbox{and all}\;\xi\in \iR^N.
\]
Suppose that $u$ is a global bounded weak solution of (\ref{GlE}).
Then $u=0$ a.e. on $\iR_+\times \iR^N$.
\end{korollar}
{\em Proof:} For $r>0$ and $x_0=0$ it follows from the proof of
Theorem \ref{Hoeldert=0} that
\begin{equation} \label{Hc2}
\omega(0,r)\le \theta \omega(0,2r),\quad r>0,
\end{equation}
where $\theta\in(0,1)$ is independent of $r$ and $u$. By induction,
(\ref{Hc2}) yields
\[
\omega(0,r)\le \theta^n \omega(0,2^nr)\le 2\theta^n
|u|_{L_\infty(\iR_+\times \iR^N)},\quad r>0,\,n\in\iN.
\]
Sending $n\to \infty$ shows that $u$ is constant. The claim then
follows by Theorem \ref{Hoeldert=0}. $\square$

$\mbox{}$

\noindent {\bf Acknowledgements:} This paper was initiated while the
author was visiting the Technical University Delft (NL) in
2003/2004. The author is greatly indebted to Philippe Cl\'ement for
many fruitful discussions and valuable suggestions.


\end{document}